\documentclass[11pt]{article}
\usepackage{latexsym, amsfonts}
\newcommand{\be}{\begin{equation}}
\newcommand{\ee}{\end{equation}}
\newcommand{\bea}{\begin{eqnarray}}
\newcommand{\eea}{\end{eqnarray}}
\newcommand{\bean}{\begin{eqnarray*}}
\newcommand{\eean}{\end{eqnarray*}}
\newcommand{\brray}{\begin{array}}
\newcommand{\erray}{\end{array}}
\newcommand{\ben}{\begin{equation}{nonumber}}
\newcommand{\een}{\end{equation}{nonumber}}

\newtheorem{dfn}{Definition}[section]
\newtheorem{thm}[dfn]{Theorem}
\newtheorem{lmma}[dfn]{Lemma}
\newtheorem{ppsn}[dfn]{Proposition}
\newtheorem{crlre}[dfn]{Corollary}
\newtheorem{xmpl}[dfn]{Example}
\newtheorem{rmrk}[dfn]{Remark}
\newcommand{\bdfn}{\begin{dfn}}
\newcommand{\bthm}{\begin{thm}}
\newcommand{\blmma}{\begin{lmma}}
\newcommand{\bppsn}{\begin{ppsn}}
\newcommand{\bcrlre}{\begin{crlre}}
\newcommand{\bxmpl}{\begin{xmpl}}
\newcommand{\brmrk}{\begin{rmrk}}
\newcommand{\edfn}{\end{dfn}}
\newcommand{\ethm}{\end{thm}}
\newcommand{\elmma}{\end{lmma}}
\newcommand{\eppsn}{\end{ppsn}}
\newcommand{\ecrlre}{\end{crlre}}
\newcommand{\exmpl}{\end{xmpl}}
\newcommand{\ermrk}{\end{rmrk}}

\newcommand{\IC}{\mathbb{C}}

\newcommand{\IN}{{I\! \! N}}

\newcommand{\IR}{\mathbb{R}}

\newcommand{\IT}{\mathbb{T}}

\newcommand{\IZ}{\mathbb{Z}}

\newcommand{\cla}{{\cal A}}
\newcommand{\clb}{{\cal B}}
\newcommand{\clc}{{\cal C}}
\newcommand{\cld}{{\cal D}}
\newcommand{\cle}{{\cal E}}
\newcommand{\clf}{{\cal F}}
\newcommand{\clg}{{\cal G}}
\newcommand{\clh}{{\cal H}}
\newcommand{\cli}{{\cal I}}

\newcommand{\clk}{{\cal K}}
\newcommand{\cll}{{\cal L}}
\newcommand{\clm}{{\cal M}}
\newcommand{\cln}{{\cal N}}

\newcommand{\clq}{{\cal Q}}

\newcommand{\cls}{{\cal S}}

\newcommand{\clu}{{\cal U}}
\newcommand{\clv}{{\cal V}}

\def\a*{{\cal A}_{h,*}}
\def\B{{\cal B}(h)}
\def\B1{{\cal B}_1(h)}
\def\b{{\cal B}^{\rm s.a.}(h)}
\def\b1{{\cal B}^{\rm s.a.}_1(h)}

\newcommand{\ot}{\otimes}

\newcommand{\raro}{\rightarrow}

\newcommand{\lgl}{\langle}
\newcommand{\rgl}{\rangle}

\def \qed {$\Box$}

\begin{document}
\[
\]
\begin{center}
{\large {\bf Quantum  Group of Orientation-preserving Riemannian Isometries }}\\
by\\
{\large Jyotishman Bhowmick {\footnote {support from the National Board of Higher Mathematics, India is gratefully acknowledged}}}\\
{and}\\
{\large Debashish Goswami {\footnote {The author gratefully acknowledges support obtained  from the Indian National Academy of Sciences through the grants for a project on `Noncommutative Geometry and Quantum Groups'.}}}\\
{\large Stat-Math Unit, Kolkata Centre,}\\
{\large Indian Statistical Institute}\\
{\large 203, B. T. Road, Kolkata 700 108, India}\\
\end{center}

\begin{abstract}
We formulate a quantum group analogue of the group of orientation-preserving Riemannian isometries of a compact Riemannian spin manifold, more generally, of a (possibly $R$-twisted  and of compact type) spectral triple. The main advantage of this formulation, which is directly in terms of the Dirac operator, is that it does not need the existence of any `good ' Laplacian as in our previous works on quantum isometry groups. Several interesting examples, including those coming from Rieffel-type deformation as well as the equivariant spectral triples on $SU_\mu(2)$ and $S^2_{\mu c}$ are dicussed.
 \end{abstract}
 
 {\bf 2000 Mathematics Subject Classification} : Primary 58B32, Secondary 16W30, 46L87, 46L89.
 
{\bf Key words and phrases}: Compact quantum group, quantum isometry groups, spectral triples
 
\section{Introduction}
Motivated by  the formulation of quantum automorphism groups by Wang (following Alain Connes' suggestion) 
(\cite{free}, \cite{wang}), and the study of their properties  by a number of mathematicians (see, e.g.
\cite{ban1}, \cite{ban2}, \cite{bichon}, \cite{univ1} and references therein), we have introduced in an earlier article (\cite{goswami}) a quantum group analogue of the group of Riemannian isometries of a classical or noncommutative manifold. In a follow-up article \cite{jyotish} we have also computed these quantum groups for a number of examples. However, our formulation of quantum isometry groups in \cite{goswami} had a major drawback from the viewpoint of noncommutative geometry since it needed a `good' Laplacian to exist.  In noncommutative geometry it is not always easy to verify such an assumption about the Laplacian, and thus it would be more appropriate to have a formulation in terms of the Dirac operator directly. This is what we aim to achieve in the present article.

   The group of Riemannian isometries of a compact Riemannian manifold $M$ can be viewed as the universal object in the
      category of all compact metrizable groups acting on $M$, with smooth and isometric action. Moreover, assume that the manifold has a spin structure (hence in particular orientable, so we can fix a choice of orientation) and $D$ denotes the conventional Dirac operator acting as an unbounded self-adjoint operator on the Hilbert space $\clh$ of square integrable spinors. Then, it can be proved that a group action on the manifold lifts as a unitary represention on the Hilbert space $\clh$ which commutes with $D$ if and only if   the action on the manifold is an orientation preserving isometric action. Therefore, to define the quantum analogue of the group of orientation-preserving Riemannian sometry group of a possibly noncommutative manifold given by a spectral triple $(\cla^\infty, \clh, D)$, 
       it is reasonable to  consider a
     category  ${\bf Q}^\prime$ of compact quantum groups having unitary (co-) representation, say $U$, on $\clh$,  which commutes with $D$, and  the action on $\clb(\clh)$ obtained by conjugation maps $\cla^\infty$ into its weak closure. 
A universal object in this category, if it exists, should define the `quantum group of orientation preserving Riemannian isometries' of the underlying spectral triple. Unfortunately, even in the finite-dimensional (but with noncommutative $\cla$) situation  this category may often fail to have a universal object, as will be discussed later. It turns out, however, that if we fix any suitable faithful functional on $\clb(\clh)$ (to be interpreted as the choice of a `volume form') then there exists a universal object in the subcategory of ${\bf Q}^\prime$ obtained by restricting the object-class to the quantum group actions which also preserve the given functional.  The  subtle point to note here is that unlike the classical group actions on $\clb(\clh)$ which always preserve the usual trace, a quantum group action may not do so. In fact, it was proved by one of the authors in \cite{goswami_rmp} that given an object $(\clq, U)$ of ${\bf Q}^\prime$ (where $\clq$ is the compact quantum group and $U$ denotes its unitary co-representation on $\clh$), we can find a  positive invertible operator $R$ in $\clh$ so that the given spectral triple is $R$-twisted in the sense of \cite{goswami_rmp} and the corresponding functional $\tau_R$ (which typically differs from the usual trace of $\clb(\clh)$ and can have a nontrivial modularity) is preserved by the action of $\clq$. This makes it quite natural to work in the setting  of twisted spectral data (as defined in \cite{goswami_rmp}).

       Motivated by the ideas of Woronowicz and Soltan, we actually consider a bigger category. The group of orientation-preserving Riemannian isometries of a classical manifold, viewed as a compact metrizable space (forgetting the group structure), can be seen to be the universal object of  a category whose object-class consists of subsets (not necessarily subgroups) of the set of such  isometries of the manifold. Then it can be proved that this universal compact set has a canonical group structure. A natural quantum analogue of this has been formulated by us, called the category of `quantum families of smooth orientation preserving Riemannian isometries'.  The underlying $C^*$-algebra of the quantum  group of orientation preserving isometries (whenever exists)  has been identified with the universal object in this bigger category and moreover, it is shown to be equipped with a canonical coproduct making it into a compact quantum group.  

We discuss a number of examples, covering both the examples coming from Rieffel-type deformation as well as the equivariant spectral triples constructed recently on $SU_\mu(2)$ and $S^2_{\mu c}$. It may be relevant to point out here that it was not clear whether one could accommodate the spectral triples on $SU_\mu(2)$ and $S^2_{\mu c}$ in the framework of our previous work on quantum isometry groups, since it is very difficult to give a nice description of the space of `noncommutative' forms and hence the Laplacian for these examples. However, the present formulation in terms of the Dirac operator makes it easy to accommodate them, and we have been able to identify $U_\mu(2)$ and $SO_\mu(3)$ as the  universal quantum group of orientation preserving isometries for the spectral triples on $SU_\mu(2)$ and $S^2_{\mu c}$ respectively (the computations for $S^2_{\mu c}$ has been presented in \cite{sphere}). 

For readers' convenience, let us briefly sketch the plan of the paper. The section 2 is devoted to the 
definition and existence of quantum group of orientation preserving isometries, which begins with a characteization of the group of such isometries in the 
 classical case (subsection 2.1), motivating the quantum formulation elaborated in subsection 2.2. The other two subsections of section 2 discuss sufficient conditions
 for ensuring a $C^*$ action of the quantum group of (orientation preserving) isometries (2.3) and for the existence of universal object without fixing a volume 
 form (2.4). Then in section 3 we study the connections of our approach with that of \cite{goswami}. Section 4 is devoted to the explicit examples and computations,
  and in the last section we skecth a general principle of computing quantum group of orientation oreserving isometries of spectral triples obtained by Rieffel
 deformation.  

We conclude this section with an important remark about the use of the phrase `orientation -preserving' in our terminology. Let us make it clear that by a `classical spectral triple' we always mean the spectral triple obtained by the Dirac operator on the spinors (so, in particular, manifolds are assumed to be compact Riemannian spin manifolds), and not just any spectral triple on the commutative algebra $C^\infty(M)$. This is absolutely crucial in view of the fact that the Hodge dirac operator $d+d^*$ on the $L^2$-space of differenttial forms also gives a spectral triple of compact type on any compact Riemannin (not necessarily with a spin structure) manifold $M$, but the the action of the full isometry group $ISO(M)$ (and not just the subgroup of orientation preserving isometries $ISO^+(M)$, even when $M$ is orientable) lifts to a  canonical unitary representation on this space commuting with $d+d^*$. In fact, the category of groups acting on $M$ such that the action comes from a unitary representation commuting with $d+d^*$,  has $ISO(M)$, and not $ISO^+(M)$, as its universal object.  So, one must stick to the Dirac operator on spinors to obtain the group of orientation preserving isometries in the usual geometric sense. This also has a natural quantum generalization, as we shall see in section 3.  

\section{Definition and existence of the quantum  group of orientation-preserving isometries}
\subsection{The classical case}
We first discuss the classical situation clearly, which will serve as a motivation for our quantum formulation. 

We begin with  a few basic facts about topologizing the space $C^\infty(M,N)$ where $M,N$ are smooth manifolds. 
Let $ \Omega $ be an open set of $\IR^n$.  We endow  $ C^{\infty} ( \Omega ) $ with the usual Frechet topology coming from  uniform convergence (over compact subsets) of  partial derivatives of all orders. The space $C^\infty(\Omega)$ is complete w.r.t. this topology, so is a Polish space in particular. Moreover, 
by the Sobolev imbedding theorem, $\cap_{k \geq 0} H_k ( \Omega ) = C^{\infty} ( {\Omega} ) $ as a set, where $H_k(\Omega)$ denotes the $k$-th Sobolev space. Thus, $C^\infty(\Omega)$ has also the Hilbertian seminorms coming from the Sobolev spaces, hence the corresponding Frechet topology. We claim that these two topologies on $C^\infty(\Omega)$ coincide.   Indeed, the inclusion map from $ C^{\infty} ( {\Omega} ) $  into $  \cap_k H_k ( \Omega )  $,  is continuous and surjective, so by the open mapping theorem for Frechet space, the inverse is also continuous, proving our claim. 

Given  two  second countable smooth manifolds $M,N$, we shall equip  $C^\infty(M,N)$ with  the weakest locally convex topology making $C^\infty(M,N) \ni \phi \mapsto f \circ \phi \in C^\infty(M)$ Frechet continuous for every $f \in C^\infty(N)$. 



For topological or smooth fibre or principal  bundles $E,F$ over a second countable smooth manifold $M$, we shall denote by ${\rm Hom}(E,F)$ the set of bundle  morphisms from $E$ to $F$.   We remark that 
 the total space of a locally trivial topological bundle such that the base and the fibre spaces are locally compact Hausdorff second countable must itself be so,  hence in particular Polish.

 In particular,  if $E$, $F$ are locally trivial principal $G$-bundles over a common base, such that the (common) base as well as structure group $G$ are locally compact Hausdorff and second countable,  then ${\rm Hom}(E,F)$ and $C(X, {\rm Hom}(E,F))$ are Polish spaces, where $X$ is a compact space. 

We need a standard fact, stated below as a lemma, about the measurable lift of Polish space valued functions. 
\blmma
\label{selection}
Let $M$ be a compact metrizable space, $B, \tilde{B}$ Polish spaces (complete separable metric spaces) such that  there is an $n$-covering map $\Lambda : \tilde{B} \raro B$.  Then any continuous map $\xi : M \raro B$ admits a lifting $\tilde{\xi} : M \raro \tilde{B}$ which is Borel measurable  and $\Lambda \circ \tilde{\xi}=\xi$.  In particular, if $\tilde{B}$ and $B$ are topological bundles over $M$, with $\Lambda$ being a bundle map, any continuous section of $B$ admits a lifting which is a measurable section of $\tilde{B}$. 
\elmma
The proof is a trivial consequence of  the selection theorem due to  Kuratowski and Ryll-Nardzewski ( see \cite{sms}, Theorem 5.2.1).

We shall now give an operator-theoretic characterization of the classical group of orientation-preserving Riemannian isometries, which will be the motivation of our definition of its quantum counterpart. Let $M$ be a compact Riemannian spin manifold, with a fixed choice of orientation.  
We note that (see, e.g. \cite{friedrich}) the spinor bundle $S$ is the associated bundle of a principal $Spin(n)$-bundle, say $P$, on $M$ ($n$=dimension of $M$), which has a canonical $2$-covering bundle-map $\Lambda$ from $ P$  to the frame-bundle $ F$ (which is an $SO(n)$-principal bundle), such that locally $\Lambda$ is of the form $({\rm id}_M \ot \lambda)$, where $\lambda : Spin(n) \raro SO(n)$ is the canonical $2$-covering group homomorphism. Let $f$ be a smooth orientation preserving Riemannian isometry of $M$, and consider the bundles $E={\rm Hom}(F, f(^*(F))$ and $\tilde{E}={\rm Hom}(P, f^*(P))$ (where ${\rm Hom}$ denotes the set of bundle maps). We view $df$ as a section of the bundle $E$ in the natural way.  By the Lemma \ref{selection} we obtain a measurable lift $\tilde{df} : M \raro \tilde{E}$, which is a measurable section of $\tilde{E}$. Using this, we define $U$ as follows. Given a (measurable) section $\xi$ of $S=P \times_{Spin(n)} \Delta_n$ (where $\Delta_n$ is as in  \cite{friedrich} ), say of the form $\xi(m)=[p(m), v]$, with $p(m) \in P_m, v \in \Delta_n$,   we define $\tilde{\xi}$ by $\tilde{\xi}(m)=[\tilde{df}(f^{-1}(m))(p(f^{-1}(m))), v] $. Note that sections of the above form constitute a total subset in $L^2(S)$, and the map $\xi \mapsto \tilde{\xi}$   is clearly a densely defined linear map on $L^2(S)$, whose  fibre-wise action is unitary since the $Spin(n)$ action is so on $\Delta_n$. Thus it extends to a unitary $U$ on $\clh=L^2(S)$.   Any such $U$, induced by the map $f$, will be denoted by $U_f$ (it is not unique since the choice of the lifting used in its construction is not unique). 
\bthm
\label{classical_case}
Let $M$ be a compact Riemannian spin manifold (hence orientable , and fix a choice of oreientation)  with the usual Dirac operator $D$ acting as  an unbounded self-adjoint operator on the Hilbert space $\clh$ of the square integrable spinors, and let $S$ denote the spinor bundle, with $\Gamma(S)$ being the $C^\infty(M)$ module of smooth sections of  $S$. Let $f: M \raro M$ be a smooth one-to-one  map which is a Riemannian isometry  preserving the orientation. Then the unitary $U_f$ on $\clh$  commutes with $D$ and  $U_f M_\phi U_f^*=M_{\phi \circ f}$, for any $\phi \in C(M)$, where $M_\phi$ denotes the operator of multiplication by $\phi$ on $L^2(S)$. Moreover, when the dimension of $ M $ is even, $ U_f $ commutes with the canonical grading $ \gamma $ on $ L^{2} ( S ).$ 

Conversely, suppose that $U$ is a unitary on $\clh$ such that $ UD = DU $ and the map $\alpha_U(X)=UXU^{-1}$ for $X \in \clb(\clh)$ maps $\cla=C(M)$ into $L^\infty(M)=\cla^{\prime \prime}$, then there is a smooth one-to-one orientation-preserving  Riemannian isometry $f$ on $M$ such that $U=U_f$.  We have the same result in the even case, if we assume furthermore that $ U \gamma = \gamma U.$


\ethm
{\it Proof:}\\
From the construction of $U_f$, it is clear that $U_f M_\phi U_f^{-1}=M_{\phi \circ f}$. Moreover, since the Dirac operator $D$ commutes  with the $Spin(n)$-action on $S$, we have $U_f D=D U_f$ on each fibre, hence on $L^2(S)$.In the even dimensional case, it is easy to see that $ Spin ( n ) $ action commutes with $ \gamma ,$ hence $ U_f $ does so.

For the converse, first note that $\alpha_U$ is a unital $\ast$-homomorphism on $L^\infty(M, dvol)$ and thus must be of the form $\psi \mapsto \psi \circ f$ for some measurable $f$. We claim that $f$ must be smooth. Fix any smooth $g$ on $M$ and consider $\phi=g \circ f$. We have to argue that $\phi$ is smooth. Let $\delta_D$ denote the generator of the strongly continuous one-parameter group of automorphism $\beta_t(X)=e^{itD} X e^{-itD}$ on $\clb(\clh)$ (w.r.t. the weak operator topology, say). From the assumption that $D$ and $U$ commute it is clear that  $\alpha_U$ maps $\cld:=\bigcap_{n \geq 1} {\rm Dom}(\delta_D^n)$ into itself, so in particular, since $C^\infty(M) \subset \cld$, we have that $\alpha_U(M_\phi)=M_{  \phi \circ g }$ belongs to $\cld$. We claim that this implies the smoothness of $\phi$. Let $m \in M$ and choose a local chart $(V, \psi)$   at $m$,  with the coordinates $(x_1,...,x_n)$, such that $\Omega=\psi(V) \in \IR^n$ has compact closure, $S|_{V}$ is trivial and $D$ has the local  expression $D= i \sum_{j=1}^n \mu(e_j) \nabla_j$, where $\nabla_j=\nabla_{\frac{\partial}{\partial x_j}}$ denotes the covariant derivative (w.r.t. the canonical Levi civita connection) operator along the vector field $\frac{\partial}{\partial x_j}$ on $L^2(\Omega)$ and $\mu(v)$ denotes the Clifford multiplication by a vector $v$. Now, $\phi \circ \psi^{-1} \in L^\infty(\Omega) \subseteq L^2(\Omega)$ and it is easy to observe from the above  local structure of $D$  that 
$[D, M_\phi]$ has the local expression $\sum_j i M_{\frac{\partial}{\partial x_j}\phi} \ot \mu(e_j)$. Thus,  
 the fact $M_\phi \in \bigcap_{n \geq 1} {\rm Dom}(\delta_D^n)$ implies  $\phi \circ \psi^{-1}  \in {\rm Dom}(d_{j_1}...d_{j_k})$ for every integer tuples $(j_1,...,j_k)$, $j_i \in \{ 1,..., n\}$, where $d_j:=\frac{\partial}{\partial x_j}$.  In other words, $\phi \circ \psi^{-1} \in H^k(\Omega)$ for all $k \geq 1$, where $H^k(\Omega)$ denotes the $k$-the Sobolev space on $\Omega$ (see \cite{rosenberg}).  By Sobolev's Theorem (see, e.g. \cite{rosenberg},Corollary 1.21, page 24)  it follows that $\phi \circ \psi^{-1} \in C^\infty(\Omega)$. 

We note that $ f $ is one-one as $ \phi \rightarrow \phi \circ f $ is an automorphism of $ L^{\infty}.$
 Now, we shall show that $f$ is an isometry of the metric space $(M,d)$, where $d$ is the metric coming from the Riemannian structure, and we have the explict formula (see \cite{con}) 
$$ d(p,q)={\rm sup}_{\phi \in C^\infty(M) , \| [D, M_\phi] \| \leq 1} |\phi(p)-\phi(q)|.$$ Since $U$ commutes with $D$, we have $\| [D, M_{\phi \circ f}]\|=\|[D, UM_\phi U^*]\|=\| U[D, M_\phi]U^* \|=\|[D, M_\phi]\|$ for every $\phi$, from which it follows that $d(f(p), f(q))=d(p,q)$.  Finally, $f$ is orientation preserving if and only if the volume form (say $\omega$) which defines the choice of orientation is preserved by the natural action of $df$, i.e. $(df \wedge .... \wedge df)(\omega)=\omega$. This will follow from the explicit description of $\omega$ in terms of $D$, given by (see \cite{Varilly_book} Page 26, also see \cite{connes_characterization} )
$$ \omega(\phi_0d\phi_1...d \phi_n)=\tau(M_{\phi_0} [D, M_{\phi_1}]...[D, M_{\phi_n}]),$$
in the odd case and  $$ \omega(\phi_0d\phi_1...d \phi_n)=\tau( \gamma M_{\phi_0} [D, M_{\phi_1}]...[D, M_{\phi_n}]),$$ in the even case where $\phi_0,...,\phi_n \in C^{\infty}(M), \gamma$ is the grading operator and $\tau$ denotes the volume integral. In fact, $\tau(X)={\rm Lim}_{t \raro 0+} \frac{{\rm Tr}(e^{-tD^2}X)}{{\rm Tr}(e^{-tD^2})}$, where ${\rm Lim}$ denotes a suitable Banach limit, which implies $\tau(UXU^*)=\tau(X)$ for all $X \in \clb(\clh)$ (using the fact that $D$ and $U$ commute). Thus, 
 
$$ \omega(\phi_0 \circ f ~d(\phi_1\circ f) \ldots d( \phi_n\circ f)) $$
$$ = \tau( \gamma U M_{\phi_{0}}U^{*} U [D, M_{\phi_{1}}]U^{*}...U[D, M_{\phi_{n}}]U^{*}) $$
$$ =  \tau(U \gamma M_{\phi_{0}} [D, M_{\phi_1}]...[D, M_{\phi_n}]U^*) $$
$$ = \tau( \gamma M_{\phi_{0}} [D, M_{\phi_1}]...[D, M_{\phi_n}]) $$
$$ = \omega(\phi_0 d\phi_1...d \phi_n).$$

with the understanding that $ \gamma = I $ when the spectral triple is odd. \qed

Now we turn to the case of a family of maps.  
We first prove a useful general fact. 

\blmma

\label{classical_case_family_ofmaps2}
Let $\cla$ be a $C^*$ algebra and $\omega, \omega_j $ ($j=1,2,...$) be states on $\cla$ such that $\omega_j \raro \omega$ in the weak-$\ast$ topology of $\cla^*$. Then for any separable Hilbert space $\clh$ and 
$ \forall ~ Y \in \clm ( \clk ( \clh ) \otimes \cla ),$ we have  $ ({\rm  id}  \otimes \omega_j ) (Y) \rightarrow ({\rm  id} \otimes \omega ) (Y) $ in the S.O.T.

\elmma

{\it Proof:}\\ Clearly,  $({\rm id} \ot \omega_j)(Y) \raro ({\rm id} \ot \omega)(Y)$ (in the strong operator topology)  for all $Y \in {\rm Fin}(\clh) \ot_{\rm alg} \cla$, where ${\rm Fin}(\clh)$ denotes the set of finite rank operators on $\clh$.  Using the strict density of ${\rm Fin}(\clh) \ot_{\rm alg} \cla $ in $\clm(\clk(\clh) \ot \cla)$, we choose, for a given $Y \in \clm(\clk(\clh) \ot \cla)$, $\xi \in \clh$ with $\| \xi \| =1$,  and $\delta>0$, an element $Y_0 \in {\rm Fin}(\clh) \ot_{\rm alg} \cla$ such that $\|(Y-Y_0) (|\xi><\xi| \ot 1) \| < \delta$. Thus, 
\bean \lefteqn{\|({\rm id} \ot \omega_j)(Y) \xi- ({\rm id} \ot \omega)(Y)\xi \|}\\
&=& \|({\rm id} \ot \omega_j)(Y (|\xi><\xi| \ot 1))\xi- ({\rm id} \ot \omega)(Y ( |\xi><\xi| \ot 1))\xi \|\\
&\leq & \|({\rm id} \ot \omega_j)(Y_0 (|\xi><\xi| \ot 1))\xi- ({\rm id} \ot \omega)(Y_0 ( |\xi><\xi| \ot 1))\xi \|\\
&+& 2  \| (Y-Y_0)(|\xi><\xi| \ot 1)\|\\
&\leq & \|({\rm id} \ot \omega_j)(Y_0 (|\xi><\xi| \ot 1))\xi- ({\rm id} \ot \omega)(Y_0 ( |\xi><\xi| \ot 1))\xi \|
+2  \delta, \eean  from which it follows that $({\rm id} \ot \omega_j)(Y) \raro ({\rm id} \ot \omega)(Y)$ in the strong operator topology. \qed

We are now ready to state and prove the operator-theoretic characterization of `set of orientation preserving isometries'. 

\bthm

\label{classical_characterization_set_orientation_preserving_isometries}

Let $ X $ be a compact metrizable space and $ \psi: X \times M \rightarrow M  $ is a map such that
 $ \psi_{x} $  defined by $ \psi_{x} ( m ) = \psi ( x, m )$  is a smooth orientation preserving Riemannian isometry and $x \mapsto \psi_x \in C^\infty(M,M)$ is continuous w.r.t. the locally convex topology of $C^\infty(M,M)$ mentioned before.

Then 
there exists a ($C(X)$-linear) unitary $U_\psi$  on the  Hilbert $C(X)$-module $\clh \ot C(X)$  such that $ \forall x \in X,~U_x:= ( id \otimes ev_x ) U_{\psi} $ is a unitary of the form $U_{\psi_x}$  on the Hilbert space $\clh$ commuting with $ D $ and $ U_x M_{\phi} U^{- 1}_{x} = U_{\phi \circ \psi^{- 1}_{x}} $. If in addition, the manifold is even dimensional, then $U_{\psi_x}$ commutes with the grading operator $ \gamma. $

Conversely, if there exists a $C(X)$-linear  unitary $ U $ on $\clh \ot C(X)$ such that $ U_x:=({\rm id} \ot {\rm ev}_x)(U) $ is a unitary commuting with $ D ~ \forall x ,$ ( and $ U_x $ commutes with the grading operator $ \gamma $ if the manifold is even dimensional )  and $ ({\rm id} \ot {\rm ev}_x) \alpha_{U} ( L^{\infty}(M) ) \subseteq L^{\infty}(M) $ for all $x \in X$,  then there exists a map $ \psi : X \times M \rightarrow M $ satisfying the conditions mentioned above such that $U=U_\psi$.

\ethm

{\it Proof:}\\ 
Consider the bundles $\hat{F}=X \times F$ and $\hat{P}=X \times P$ over $X \times M$, with fibres at $(x,m)$ isomorphic with (respectively) $P_m$ and $F_m$, and where  $ F $ and $ P $ are as in Theorem \ref{classical_case}. Moreover, denote by $\Psi$ the map from $X \times M$ to itself given by $(x,m) \mapsto (x, \psi(x,m))$, and let $B=C(X, {\rm Hom}(\hat{F}, \Psi^*(\hat{F})))$, $\tilde{B}:=C(X, {\rm Hom}(\hat{P},\Psi^*(\hat{P})))$. The covering map from $P$ to $M$ induces a covering map from $\tilde{B}$ to $B$ as well.  Let $ d^{'}_\psi: M  \raro B$ be  the map given by  $ d^{'}_\psi ( m)|_{(x,m) }= d \psi_{x}|_{m}.$ Then by Theorem \ref{selection} there exists a measurable lift of $ d^{'}_\psi,$ say $ \widetilde{d^{'}_\psi} $ from $M$ into $\tilde{B} $. Since $d^{'}_\psi ( m)|_{(x,m) } \in {\rm Hom}(F_m, F_{\psi(x,m)})$, it is clear that the lift $\widetilde{d^\prime_\psi}(m)|_{(x,m)}$ will be an element of ${\rm Hom}(P_m, P_{\psi(x,m)})$.

We can identify  $  \clh \otimes C ( X )  $ with $ C ( X \rightarrow  \clh ) ,$ and since $\clh$ has a total set $\clf$ (say)  consisting of sections of the form $[p(\cdot), v]$, where $p: M \raro P$ is a measurable section of $P$ and $v \in \Delta_n$, we have a total set $\tilde{\clf}$ of $ \clh \ot C(X)$ consisting of $\clf$ valued continuous functions from $X$. Any such function can be written as   $[\Xi,v]$ with $\Xi : X \times M  \raro P$, $v \in \Delta_n$, and $\Xi(x, m) \in P_m$,  and  we define $U$ on $\tilde{\clf}$ by $U[\xi, v]=[\Theta, v]$, $$   \Theta( x,m ) =  \widetilde{d^{'}_\psi}( m )|_{ ( x,\psi^{- 1}_{x} ( m ) ) }( \Xi (x, \psi^{- 1}_{x} ( m ) ) ).$$ It is clear from  the construction of the lift that $U$ is indeed a $C(X)$-linear isometry which maps the total set $\tilde{\clf}$ onto itself, so extends to a unitary on the whole of $ \clh \ot C(X)$ with the desired properties.


Conversely, given $U$ as in the statement of the converse part of the theorem, we observe that for each $x \in X$, by Theorem \ref{classical_case},  $ ( id \otimes ev_x ) U = U_{\psi_{x}} $ for some $ \psi_x $ such that $ \psi_x $ is a smooth  orientation preserving Riemannian isometry. This defines the map $\psi$ by setting $\psi(x,m)=\psi_x(m)$. The proof will be complete if we can  show that $x \mapsto \psi_x \in C^\infty(M,M)$ is continuous, which is equivalent to showing that whenever $x_n \raro x$ in the topology of $X$, we must have $ \phi \circ \psi_{x_n} \rightarrow \phi \circ \psi_x $ in the Frechet topology of $ C^{\infty}(M)$,  for any $\phi \in C^\infty(M)$. However, by  Lemma \ref{classical_case_family_ofmaps2}, we have $ ( id \otimes ev_{x_n} ) \alpha_U ( [ D, ~M_\phi ] ) \rightarrow ( id \otimes ev_x ) \alpha_U ( [ D, ~ M_\phi ] ) $ in the S.O.T. Since  $ U $ commutes with $ D ,$ this implies $$ ( id \otimes ev_{x_n} ) [ D \otimes id, ~ \alpha_{U} ( M_\phi ) ] \rightarrow ( id \otimes ev_x ) [ D \otimes id, ~ \alpha_{U} ( M_\phi ) ] ,$$ i.e. 
$$[ D, M_{\phi \circ \psi_{x_n}} ] \xi \rightarrow [ D, ~ M_{\phi \circ \psi_x} ] \xi ~{\rm  in}~  L^2 ~  \forall ~ \xi \in L^2 ( S ) .$$ 
By choosing $\phi$ which has support in a local trivializing coordinate neighborhood for $S$, and then using the local expression of $D$ as in the proof of Theorem \ref{classical_case},
we conclude that 
$ d_k ( \phi \circ \psi_{x_n} ) \rightarrow  d_k( \phi \circ \psi_x ) $ in $ L^2 $ (where $d_k$ is as in the proof of Theorem \ref{classical_case}).
Similarly, by taking repeated commutators with $D$, we can show the convergence with $d_k$ replaced by $d_{k_1}...d_{k_m}$ for any finite tuple $(k_1,...,k_m)$.  In other words, $\phi \circ \psi_{x_n} \raro \phi \circ \psi_x$ in the topology of $C^\infty(M)$ described before. \qed

\subsection{Quantum group of orientation-preserving isometries of an $R$-twisted spectral triple}
 We begin by   recalling the definition of compact quantum groups and their actions from   \cite{woro}, \cite{woro1}.  A 
compact quantum group (to be abbreviated as CQG from now on)  is given by a pair $(\cls, \Delta)$, where $\cls$ is a unital separable $C^*$ algebra 
equipped
 with a unital $C^*$-homomorphism $\Delta : \cls \raro \cls \otimes \cls$ (where $\otimes$ denotes the injective tensor product)
  satisfying \\
  (ai) $(\Delta \ot id) \circ \Delta=(id \ot \Delta) \circ \Delta$ (co-associativity), and \\
  (aii) the linear span of $ \Delta(\cls)(\cls \ot 1)$ and $\Delta(\cls)(1 \ot \cls)$ are norm-dense in $\cls \ot \cls$. \\
  It is well-known (see \cite{woro}, \cite{woro1}) that there is a canonical dense $\ast$-subalgebra $\cls_0$ of $\cls$, consisting of the matrix coefficients of
   the finite dimensional unitary (co)-representations (to be defined shortly) of $\cls$, and maps $\epsilon : \cls_0 \raro \IC$ (co-unit) and
   $\kappa : \cls_0 \raro \cls_0$ (antipode)  defined
    on $\cls_0$ which make $\cls_0$ a Hopf $\ast$-algebra.

    We say that  the compact quantum group $(\cls,\Delta)$ (co)-acts on a unital $C^*$ algebra $\clb$,
    if there is a  unital $C^*$-homomorphism (called an action) $\alpha : \clb \raro \clb \ot \cls$ satisfying the following :\\
    (bi) $(\alpha \ot id) \circ \alpha=(id \ot \Delta) \circ \alpha$, and \\
    (bii) the linear span of $\alpha(\clb)(1 \ot \cls)$ is norm-dense in $\clb \ot \cls$.\\
     \vspace{1mm}\\
It is known ( see, for example,  \cite{wang} , \cite{podles_subgroup} ) that  (bii) is equivalent to the existence of a norm-dense, unital $\ast$-subalgebra $\clb_0$ of $\clb$ such that $\alpha(\clb_0) \subseteq \clb_0 \ot_{\rm alg} \cls_0$ and on $\clb_0$, $({\rm id }\ otimes \epsilon) \circ \alpha={\rm id}$.

We shall sometimes say that $ \alpha $ is a  'topological' or $ C^* $ action to distinguish it from a normal action of Von Neumann algebraic quantum group.
      
  \bdfn
  
   A unitary ( co ) representation of a compact quantum group $ ( S, \Delta ) $ on a Hilbert space $ \clh $ is a map $ U $ from $ \clh $ to the Hilbert $\cls$ module $ \clh \otimes \cls $  such that the  element $ \widetilde{U} \in \clm ( \clk ( \clh ) \otimes \cls ) $ given by $\widetilde{U}( \xi \ot b)=U(\xi)(1 \ot b)$ ($\xi \in \clh, b \in \cls $) is a unitary satisfying  $$ ({\rm  id} \otimes \Delta ) \widetilde{U} = {\widetilde{U}}_{(12)} {\widetilde{U}}_{(13)},$$ where for an operator $X \in \clb(\clh_1 \ot \clh_2)$ we have denoted by $X_({12})$ and $X_({13})$ the operators $X \ot I_{\clh_2} \in \clb(\clh_1 \ot \clh_2 \ot \clh_2)$, and $\Sigma_{23} X_({12}) \Sigma_{23}$ respectively ($\Sigma_{23}$ being the unitary on $\clh_1 \ot \clh_2 \ot \clh_2$ which flips the two copies of $\clh_2$).
  
Given a unitary representation $U$ we shall denote by $\alpha_U$ the $\ast$-homomorphism $\alpha_U(X)=\widetilde{U}(X \ot 1){\widetilde{U}}^*$ for $X \in \clb(\clh)$. For a  not necessarily bounded, densely defined (in the weak operator topology)  linear functional $\tau$ on $\clb(\clh)$,  we say that $\alpha_U$ preserves $\tau$ if $\alpha_U$ maps a suitable (weakly) dense $\ast$-subalgebra   (say $\cld$) in the domain of $\tau$ into $\cld \ot_{\rm alg} \cls$ and $( \tau \ot {\rm id}) (\alpha_U(a))=\tau(a)1_\cls$  for all $a \in \cld$. When $\tau$ is bounded and normal, this is equivalent to $(\tau \ot {\rm id}) (\alpha_U(a))=\tau(a) 1_\cls$ for all $a \in \clb(\clh)$. 

We say that a (possibly unbounded) operator $T$ on $\clh$ commutes with $U$ if $T \ot I$ (with the natural domain) commutes with $\widetilde{U}$. Sometimes such an operator will be called $U$-equivariant.
\edfn

 Let us now recall the concept of universal quantum groups as in
\cite{univ1}, \cite{free}
and references therein. We shall use most of the terminologies of
\cite{free}, e.g. Woronowicz $C^*$ -subalgebra, Woronowicz
$C^*$-ideal etc, however with the exception that we shall call the
Woronowicz $C^*$ algebras just compact quantum groups, and not use
the term compact quantum groups for the dual objects as done in
\cite{free}.  
For an $n \times n$ positive invertible matrix $Q=(Q_{ij})$. let $A_u(Q)$ be the compact quantum group defined and  studied in \cite{wang}, \cite{univ1}, which is the universal $C^{*}$-algebra generated by $ \{ u^{Q}_{kj}, k,j=1,...,d_{i} \}$ such that $u:=(( u^{Q}_{kj} ))$ satisfies \be \label{wangalg} u u^*=I_n =u^{*}u, ~~u^{\prime} Q  \overline{u} Q^{-1}=I_n=Q{\overline{u}} Q^{-1} u^{\prime}.\ee Here 
$u^{\prime} =(( u_{ji} ))$ and $\overline{u}=(( u_{ij}^{*} ))$. The coproduct, say $\tilde{\Delta}$, is given by, $$ \tilde{\Delta}(u_{ij})=\sum_k u_{ik} \ot u_{kj}.$$ It may be noted that $A_u(Q)$ is the universal object in the category of compact quantum groups which admit an action on the finite dimensional $C^*$ algebra $M_n(\IC)$ which preserves  the functional $M_n \ni x \mapsto {\rm Tr({Q}^{T} x)}$,( see \cite{wangergodic} ) where we
refer the reader to \cite{univ1}  for a detailed  discussion on the structure and classification of
such quantum groups.

In view of the  characterization of orientation-preserving isometric action on  a classical manifold ( Theorem \ref{classical_characterization_set_orientation_preserving_isometries} ), we give the following definitions.
         \bdfn 
         \label{def_q_fam}A quantum family of orientation preserving  isometries for the ( odd ) spectral triple $({\cla^\infty}, \clh, D)$ is given by a pair $(\cls, U)$ where $\cls$ is a separable unital $C^*$-algebra and  $U$ is a linear map from $\clh$ to $\clh \ot \cls$ such that $\widetilde{U}$ given by $\widetilde{U}( \xi \ot b)=U(\xi) (1 \ot b)$ $(\xi \in \clh$, $b \in \cls$) extends to a unitary element of  $ \clm(\clk(\clh) \ot \cls)$ satisfying the following \\
(i) for every state $\phi$ on $\cls$ we have $U_\phi D=DU_\phi$, wher $U_\phi:=({\rm id} \ot \phi)(\widetilde{U})$;\\
(ii) $({\rm id} \ot \phi) \circ \alpha_U(a) \in ({\cla^\infty})^{\prime \prime}$ $\forall a \in \cla^\infty$ for every state $\phi$ on $\cls$, where $\alpha_U(x):=\widetilde{U}( x \ot 1) {\widetilde{U}}^* $ for $x \in \clb(\clh)$.

In case the $C^*$-algebra $\cls$ has a coproduct $\Delta$ such that $(\cls,\Delta)$ is a compact quantum group and $U$ is a unitary representation  of $(\cls, \Delta)$ on $\clh$, we say that $(\cls, \Delta)$ acts  by orientation-preserving isometries  on the spectral triple.

In case the spectral triple is even with grading operarator $ \gamma ,$ a quantum family of orientation preserving isometries $ ( \cla^{\infty}, \clh, D, \gamma ) $ will be defined exactly as above,  with the only extra condition being that $ U $ commutes with $ \gamma. $

\edfn

From now on, we will mostly consider odd spectral triples. However let us remark that in the even case,  all the definitions and results obtained by us will go through with some obvious modifications.

Consider the category ${\bf Q}\equiv {\bf Q}(\cla^\infty, \clh, D)\equiv {\bf Q}(D)$ with the object-class consisting of all quantum families of 
orientation preserving isometries $(\cls, U)$ of the given spectral triple, and the set of morphisms ${\rm Mor}((\cls,U),(\cls^\prime,U^\prime))$ 
being the set of unital $C^*$-homomorphisms $\Phi : \cls \raro \cls^\prime$ satisfying $({\rm id} \ot \Phi) (U)=U^\prime$. 
We also consider another category ${\bf Q}^\prime \equiv {\bf Q}^\prime(\cla^\infty, \clh, D) \equiv {\bf Q}^\prime(D)$ whose objects are triplets $(\cls, \Delta, U)$,  where $(\cls,\Delta)$ is a compact quantum group acting by orientation preserving isometries on the given spectral triple, with $U$ being the corresponding unitary representation. The morphisms  are the homomorphisms of compact quantum groups which are also morphisms of the underlying quantum families of orientation preserving isometries. The forgetful functor $F: {\bf Q}^\prime \raro {\bf Q}$ is clearly faithful, and we can view $F({\bf Q}^\prime)$ as a subcategory of ${\bf Q}$.

It is easy to see that any object $(\cls, U)$  of the category ${\bf Q}^\prime$ gives an equivariant spectral triple $(\cla^\infty, \clh, D)$  w.r.t. the action of $\cls$ implemented by $U$. It  may be noted that recently there has a lot of interest and work (see, for example, \cite{partha}, \cite{con2}, \cite{landiqgp}) towards construction of quantum group equivariant spectral triples. In all these works, given a $C^*$-algebra $\cla \subseteq \clb(\clh)$ and a CQG $\clq$ having a unitary representation $U$ on $\clh$ such that ${\rm ad}_U$ gives an action of $\clq$ on $\cla$,  the authors investigate the possibility of constructing a (nontrivial) spectral triple $(\cla^\infty, \clh, D)$ on a suitable dense subalgebra $\cla^\infty$ of   $\cla$ such that $\tilde{U}$ commutes with $D \ot 1$, i.e. $D$ is equivariant. Our interest here is  in the (sort of) 
 converse direction: given a spectral triple, we want to consider all possible CQG represenations which w.r.t. which the spectral triple is equivariant; and if there exists a universal object in the corresponding category, i.e ${\bf Q}^\prime$, we should call it the quantum group of orientation preserving isometries. 

 Unfortunately, in general ${\bf Q}^\prime$ or ${\bf Q}$ will not have a universal object. It is easily seen by taking the standard example $\cla^\infty=M_n(\IC)$, $\clh=\IC^n$, $D=I$. However, the fact that comes to our rescue is that a universal object exists in each of the subcategories which correspond to the CQG actions preserving a given faithful functional on $M_n$.
 
  On the other hand, given any equivariant spectral triple, it has been shown in \cite{goswami_rmp} that there is a (not necessarily unique) canonical faithful functional which is preserved by the CQG action.  For readers' convenience, we state this result  (in a form suitable to us) briefly here:
\bppsn
\label{5678}
Given a spectral triple $(\cla^\infty, \clh, D)$ (of compact type)  which is $\clq$-equivariant w.r.t. a representation of a CQG $\clq$ on $\clh$, we can construct a positive (possibly unbounded) invertible operator $R$ on $\clh$ such that $(\cla^\infty, \clh, D, R)$ is  a twisted spectral data, i.e.  \\
(a) $R$ commutes with $D$ and $\forall s \in \IR$, \\
(b) the map $a \mapsto \sigma_s(a):=R^{-s} a R^s$ gives an automorphism of $\cla^\infty$ (not $\ast$-preserving) satisfying  $\sup_{s \in [-n,n]} \| \sigma_s(a) \| < \infty$ for all positive integer $n$;\\
 and moreover, we have \\
(c) $\alpha_U$ preserves the functional $\tau_R$ (defined at least on a weakly dense $\ast$-subalgebra $\cle_D$ of $\clb(\clh)$ generated by the rank-one operators of the form $|\xi><\eta|$ where $\xi, \eta$ are eigenvectors of $D$) given  by $$ \tau_R(x)= Tr ( R x ),~~x \in \cle_D.$$  
\eppsn

\brmrk
If $V_\lambda$ denotes the eigenspace  of $D$ corresponding to the eigenvalue, say $ \lambda ,$  it is clear that $ \tau_{R} (X) = e^{t \lambda^2} {\rm Tr}(Re^{-tD^2}X) $ for all $X =|\xi><\eta|$ with $\xi, \eta \in V_\lambda$ and for any $t>0$. Thus, the $\alpha_U$-invariance of the functional $\tau_R$ on $\cle_D$ is equivalent to the $\alpha_U$-invariance of the functional $X \mapsto {\rm Tr}(X Re^{-tD^2})$ on $\cle_D$ for each $t>0$. 
 If, furthermore, the $R$-twisted spectral triple is $\Theta$-summable in the sense that  $Re^{-tD^2}$ is trace class for every $t>0$, the above is also equivalent to the $\alpha_U$-invariance of the bounded normal functional $X \mapsto {\rm Tr}(X Re^{-tD^2})$ on the whole of $\clb(\clh)$.  In particular, this implies that  $\alpha_U$ preserves the  functional    $ \clb(\clh) \ni x \mapsto {\rm Lim}_{t \raro 0+} \frac{{\rm Tr}(xRe^{-tD^2})}{{\rm Tr}(Re^{-tD^2})},$ where ${\rm Lim}$ is a suitable Banach limit discussed in, e.g \cite{fro}.  
\ermrk

This motivates the following definition:
\bdfn
Given an  $R$-twisted spectral triple $(\cla^{\infty}, \clh, D, R)$, a quantum family  of orientation preserving isometries $(\cls, U)$ of $(\cla^\infty, \clh, D)$ is said to preserve the $R$-twisted volume, (simply said to be  volume-preserving if $R$ is understood) if one has  $(\tau_R \otimes  {\rm id} ) (\alpha_U(x))= \tau_R(x)1_\cls$ for all $x \in \cle_D$, where $\cle_D$ and $\tau_R$ are as in Proposition \ref{5678}. We shall also call $(\cls, U)$ a quantum family of orientation-preserving isometries of the $R$-twisted spectral triple. 

If, furthermore, the  $C^*$-algebra $\cls$ has a coproduct $\Delta$ such that $(\cls,\Delta)$ is a CQG and $U$ is a unitary representation  of $(\cls, \Delta)$ on $\clh$, we say that $(\cls, \Delta)$ acts  by volume and orientation-preserving isometries  on the $R$-twisted spectral triple.

 We shall consider the categories ${\bf Q}_R \equiv {\bf Q}_R(D)$ and ${\bf Q}^\prime_R \equiv {\bf Q}^\prime_R(D)$ which are the full subcategories of ${\bf Q}$ and ${\bf Q}^\prime$ respectively, obtained by restricting the object-classes to the volume-preserving quantum families. 
\edfn

\brmrk
We shall not need the full strength of the definition of twisted spectral data here; in particular the condition (b) in the Proposition \ref{5678}. However, we shall continue to work with the usual definition of $R$-twisted spectral data, keeping in mind that all our results are valid without assuming ( b ).
\ermrk

Let us now fix a spectral triple $(\cla^\infty, \clh, D)$ which is of compact type. The $C^*$-algebra generated by $\cla^\infty$ in $\clb(\clh)$ will be denoted by $\cla$. Let $\lambda_0=0, \lambda_1, \lambda_2,\cdot \cdot \cdot  $ be the eigenvalues of $D$ with $V_i$ denoting the ($d_i$-dimensional, $d_i<\infty$) eigenspace for $\lambda_i$. Let $\{ e_{ij}, j=1,..., d_i \}$ be an orthonormal basis of $V_i$. 
We also assume that there is a positive  invertible $R$ on $\clh$ such that $(\cla^\infty, \clh, D, R)$ is 
 an $R$-twisted   spectral triple. $R$ must have the form $R|_{V_i}=R_i, say$, with $R_i$ positive invertible $d_i \times d_i$ matrix. Let us denote the CQG $A_u(R^T_i)$ by $\clu_i$, with its canonical unitary representation $\beta_i$ on $V_i \cong \IC^{d_i}$, given by $\beta_i(e_{ij})=\sum_k e_{ik} \ot u^{R^T_i}_{kj}$. Let $\clu$ be the free product of $\clu_i, i=1,2,...$ and $\beta=\ast_i \beta_i$ be the corresponding free product representation of $\clu$ on $\clh$. We shall also consider the  corresponding  unitary element $\tilde{\beta}$ in $\clm(\clk(\clh) \ot \clu)$.
 \blmma
\label{lem2} Consider the $R$-twisted spectral triple $(\cla^\infty,\clh,D)$  and let $(\cls,U)$ be a quantum family of volume and orientation preserving  isometries of the given spectral triple.  Moreover, assume that the map  $U$ is faithful in the sense that there is no
proper  $C^*$-subalgebra $\cls_1$ of $\cls$ such that
$\widetilde{U} \in \clm(\clk(\clh ) \ot \cls_1)$. 
Then  we
can find a $C^*$-isomorphism  $\phi : \clu/
\cli \raro \cls$ between $\cls$ and a quotient of $\clu$ by a
 $C^*$-ideal
 $\cli$ of $\clu$, such that $ U= ({\rm id}\ot \phi) \circ ({\rm id} \ot \Pi_\cli) \circ
 \beta$, where $\Pi_\cli$ denotes the quotient map from $\clu$ to
 $\clu/\cli$.

 If, furthermore, there is a compact quantum group structure on $\cls$ given by a coproduct  $\Delta$ such that $(\cls,\Delta, U)$ is an object in ${\bf Q}^\prime_R$, the ideal $\cli$ is a Woronowicz $C^*$-ideal and the $C^*$-isomorphism $\phi : \clu/ \cli \raro \cls$ is a morphism of compact quantum groups.
\elmma
{\it Proof :}\\

It is clear that $U$  maps $V_i $ into $V_i \otimes
\cls$ for each $i$. Let $v^{(i)}_{kj}$ ($j,k=1,...,d_i$) be the
elements of $\cls$ such that $U(e_{ij})=\sum_k e_{ik} \otimes
v^{(i)}_{kj}$.  Note that $v_i:=((v^{(i)}_{kj} ))$ is a unitary in
$M_{d_i}(\IC) \otimes \cls$. Moreover, the $\ast$-subalgebra
generated by all $ \{ v^{(i)}_{kj}, i \geq 0, ,j,k \geq 1 \}$ must be
dense in $\cls$ by the assumption of faithfulness.

Consider the $\ast$-homomorphism $\alpha_i $ from the finite dimensional $C^*$ algebra $\cla_i \cong M_{d_i}(\IC)$ generated by the rank one operators $\{ |e_{ij}><e_{ik}|, j,k=1,..., d_i \}$ to $\cla_i \ot \cls$ given by $\alpha_i(y)=\widetilde{U}(y \ot 1){\widetilde{U}}^*|_{V_i}$. Clearly,  the restriction of the  functional $\tau_R$ on $\cla_i$ is nothing but the  functional given by ${\rm Tr}( R_i~ \cdot)$, where ${\rm Tr}$ denotes the usual trace of matrices. Since $\alpha_i$ preserves this functional by assumption, we get,  by the universality of $\clu_i$,  a
$C^*$-homomorphism from $\clu_i$ to $\cls$ sending $u^{(i)}_{kj} \equiv u_{kj}^{R^T_i} $
to $v^{(i)}_{kj}$, and by definition of the free product, this
induces a $C^*$-homomorphism, say $\Pi$, from $\clu$ onto $\cls$,
 so that $\clu/\cli
\cong \cls$, where  $\cli:={\rm Ker}(\Pi)$.

In case $\cls$ has a coproduct $\Delta$ making it into a compact quantum group and $U$  is a quantum group representation, it is easy to see that  the subalgebra of $\cls$ generated by $\{ v^{(i)}_{kj},~i \geq 0, j,k=1,...,d_i \}$ is a Hopf
algebra, with $\Delta(v^{(i)}_{kj})=\sum_l v^{(i)}_{kl} \ot v^{(i)}_{lj}$. From this, it follows that  $\Pi$ is Hopf-algebra morphism, hence $\cli$ is a Woronowicz $C^*$-ideal.
\qed \vspace{1mm}\\

Before we state and prove the main theorem, let us note the following elementary fact about $C^*$-algebras.

  \blmma \label{lim} Let $\clc$ be a unital 
$C^*$ algebra and $\clf$ be a nonempty collection of $C^*$-ideals (closed two-sided ideals) of
$\clc$. Let  $\cli_0$ denote the intersection of all $\cli$ in $\clf$, and  let $\rho_\cli$ denote the map $\clc /\cli_0 \ni x+\cli_0 \mapsto x+I \\ \in \clc/\cli$ for $\cli \in \clf$.  Denote by $\Omega$  the set  $\{ \omega \circ \rho_\cli, \cli \in \clf,~\omega ~{\rm state}~{\rm on}~ \clc/\cli \}$, and let $K$ be the weak-$\ast$ closure of the convex hull of $\Omega \bigcup (- \Omega)$. Then $K$ coincides with the set of bounded linear functionals $\omega$ on $\clc/\cli_0$ satisfying $\| \omega \|=1$ and $\omega(x^*+ \cli_0)=\overline{\omega(x+ \cli_0)}$.  \elmma
 {\it Proof :}\\
We have, by Lemma 4.6 of \cite{goswami} that for any $x \in \clc$, 
$$ \sup_{\cli \in \clf}   \| x+\cli \| = \|x+ \cli_0 \|,$$
where $\|
x+\cli \|=inf \{ \| x -y \|~:~y \in \cli \}$ denotes the norm in
$\clc/\cli$. Clearly, $K$ is a weak-$\ast$ compact, convex subset of the unit ball $(\clc/\cli_0)^*_1$ of the dual of $\clc/\cli_0$, satisfying $-K=K$.  If $K$ is strictly smaller than the self-adjoint part of unit ball of the dual of $\clc/\cli_0$, we can find a state $\omega$ on $\clc/\cli_0$ which is not in $K$. Considering the real Banach space $X=(\clc/\cli_0)^*_{\rm s.a.}$ and using standard separation theorems for real Banach spaces (e.g. Theorem 3.4 of \cite{rudin}, page 58), we can find a self-adjoint element $x$ of $\clc$ such that $\| x+\cli_0 \|=1,$ and   $$ \sup_{\omega^\prime \in K} \omega^\prime(x+\cli_0) < \omega(x+\cli_0).$$ Let $\gamma \in \IR$ be such that $\sup_{\omega^\prime \in K} \omega^\prime(x+\cli_0)< \gamma  < \omega(x+\cli_0).$ Fix  $0<\epsilon<\omega(x+\cli_0)-\gamma$, and let $\cli \in \clf$ be such that $\| x+\cli_0\| -\frac{\epsilon}{2} \leq \| x+\cli \| \leq \| x+\cli_0\|.$ Let $\phi$ be a state on $\clc/\cli$ such that $\| x+ \cli \|=| \phi(x+\cli)|.$ Since $x$ is self-adjoint, either $\phi(x+\cli)$ or $-\phi(x+\cli)$ equals $\| x+\cli \|$, and $\phi^\prime:=\pm \phi \circ \rho_\cli $, where the sign is chosen so that $\phi^\prime(x+\cli_0)=\| x+\cli\|$. Thus, $\phi^\prime \in K$, so $\| x+\cli_0\|=\phi^\prime(x+ \cli) \leq \gamma < \omega(x+\cli_0)-\epsilon.$ But this implies $\|x+\cli_0\| \leq \| x+\cli \| +\frac{\epsilon}{2}<\omega(x+\cli_0)-\frac{\epsilon}{2}\leq \| x+\cli_0\| -\epsilon$ (since $\omega$ is a state), which is a contradiction completing the proof.
 \qed
 \bthm
\label{main} For any  $R$-twisted  spectral triple $(\cla^\infty, \clh, D)$, the category ${\bf Q}_R$ of quantum families of volume and orientation preserving isometries has a universal (initial) object, say  $(\widetilde{\clg}, U_0)$. Moreover, $\widetilde{\clg}$ has a coproduct $\Delta_0$ such that $(\widetilde{\clg},\Delta_0)$ is a compact quantum group and $(\widetilde{\clg},\Delta_0,U_0)$ is a universal object in the category ${\bf Q}^\prime_0$.  The representation  $U_0$ is faithful.
 \ethm
 {\it Proof :}\\
Recall the $C^*$-algebra $\clu$ considered before, along with  the map $\beta $ and the corresponding  unitary  $\widetilde{\beta} \in \clm(\clk(\clh) \ot \clu ) $. For any $C^*$-ideal $\cli$ of $\clu$, 
we shall denote by $\Pi_\cli$ the canonical quotient map from $\clu$ onto $\clu/\cli$, and let $\Gamma_\cli=({\rm id} \ot \Pi_\cli) \circ \beta$. Clearly,   $\widetilde{\Gamma_\cli}=({\rm id} \ot \pi_\cli)\circ \widetilde{\beta}$ is a unitary element of  $\clm(\clk(\clh) \ot \clu/\cli )$. 
Let $\clf$ be the collection of all those  $C^*$-ideals
$\cli$ of $\clu$ such that  $({\rm id} \ot \omega) \circ \alpha_{\Gamma_\cli } \equiv ({\rm id} \ot \omega) \circ {\rm ad}_{\widetilde{\Gamma_\cli} }$ maps $\cla^\infty$ into $\cla^{\prime \prime}$ for every state $\omega$ (equivalently, every bounded linear functional) on $\clu/\cli$.  This
collection is nonempty, since the trivial one-dimensional $C^*$-algebra $\IC$ gives an object in ${\bf Q}_R$ and by Lemma \ref{lem2} we do get a member of $\clf$. Now,
let $\cli_0$ be the intersection of all ideals in $\clf$. We claim
that $\cli_0$ is again a member of $\clf$.  
Indeed, in the notation of Lemma \ref{lim}, it is clear that for $a \in \cla^\infty$, $({\rm id} \ot \phi)\circ \widetilde{\Gamma}_{\cli_0}(a) \in \cla^{\prime \prime}$ for all $\phi $ in the convex hull of $ \Omega \bigcup (-\Omega)$. Now, for any state $\omega$ on $\clu/\cli_0$, we can find, by Lemma \ref{lim}, a net $\omega_j $ in the above convex hull  (so in particular $\| \omega_j \| \leq 1~\forall j$) such that $\omega(x+\cli_0)=\lim_j \omega_j(x+ \cli_0)$ for all $x \in \clu/\cli_0$. 

 It follows from  Lemma \ref{classical_case_family_ofmaps2} that $({\rm id} \ot \omega_j)(X) \raro ({\rm id} \ot \omega)(X)$ (in the strong operator topology)  for all $X \in \clm(\clk(\clh) \ot \clu/\cli_0)$. 
Thus, for $a \in \cla^\infty$, $ ({\rm id} \ot \omega)\circ {\rm ad}_{ \widetilde{\Gamma}_{\cli_0}}(a)$ is the S.O.T. limit of $({\rm id} \ot \omega_i)\circ {\rm ad}_{ \widetilde{\Gamma}_{\cli_0}}(a)$, hence belongs to $\cla^{\prime \prime}$.

We now show   that $(\widetilde{\clg}:=\clu/\cli_0, \Gamma_{\cli_0})$ is a 
universal object in ${\bf Q}_R$. To see this, consider any object $(\cls, U)$ of ${\bf Q}_R$.  Without loss
of generality we can assume $U$ to be faithful, since
otherwise we can replace $\cls$ by the  $C^*$-subalgebra
generated by the elements $\{ v^{(i)}_{kj} \}$ appearing in the proof of Lemma \ref{lem2}.  But by
Lemma \ref{lem2} we can further assume that $\cls$ is isomorphic
with $\clu/\cli$ for some $\cli \in \clf$. Since $\cli_0 \subseteq
\cli$, we have a $C^*$-homomorphism from
$\clu/\cli_0$ onto $\clu/\cli$, sending $x+\cli_0$ to $x +\cli$, which is clearly a morphism in the category ${\bf Q}_R$. This is indeed the unique such morphism, since it is uniquely determined on the dense subalgebra generated by $\{ u^{(i)}_{kj}+ \cli_0,~i\geq 0,~j,k \geq 1 \}$ of $\widetilde{\clg}$.

To construct the coproduct on $\widetilde{\clg}=\clu/\cli_0$, we first consider $U^{(2)}: \clh \raro \clh \ot \widetilde{\clg} \ot \widetilde{\clg}$ given by   $$U^{(2)}=(\Gamma_{\cli_0})_{(12)}(\Gamma_{\cli_0})_{(13)},$$ where $U_{ij}$ is the usual `leg-numbering notation'. It is easy to see that $(\widetilde{\clg} \ot \widetilde{\clg}, U^{(2)})$ is an object in the category ${\bf Q}_R$, so by the universality of $(\widetilde{\clg},\Gamma_{\cli_0})$, we have a unique unital $C^*$-homomorphism $\Delta_0 : \widetilde{\clg} \raro \widetilde{\clg} \ot \widetilde{\clg}$ satisfying $$ ({\rm id} \ot \Delta_0) (\Gamma_{\cli_0}) =U^{(2)}.$$ 
Letting both sides act on  $e_{ij}$, we get $$ \sum_l e_{il} \ot (\pi_{\cli_0} \ot \pi_{\cli_0}) \left( \sum_k u^{(i)}_{lk} \ot u^{(i)}_{kj} \right)=\sum_l e_{il} \ot \Delta_0(\pi_{\cli_0}(u^{(i)}_{lj})).$$ Comparing coefficients of $e_{il}$, and recalling that $\tilde{\Delta}(u^{(i)}_{lj})=\sum_k u^{(i)}_{lk} \ot u^{(i)}_{kj}$ (where $\tilde{\Delta}$ denotes the coproduct on $\clu$), we have \be \label{coprod}(\pi_{\cli_0} \ot \pi_{\cli_0}) \circ \tilde{\Delta}=\Delta_0 \circ \pi_{\cli_0} \ee on the linear span of $\{ u^{(i)}_{jk}, i \geq 0,~j,k \geq 1 \}$, and hence on the whole of $\clu$. 
This implies that $\Delta_0$ maps $\cli_0={\rm Ker}(\pi_{\cli_0})$ into ${\rm Ker}(\pi_{\cli_0} \ot \pi_{\cli_0})=(\cli_0 \ot 1 + 1 \ot \cli_0) \subset \clu \ot \clu$. In other words, $\cli_0$ is a Hopf $C^*$-ideal, and hence $\widetilde{\clg}=\clu/\cli_0$ has the canonical compact quantum group structure as a quantum subgroup of $\clu$. It is clear from the relation (\ref{coprod}) that  $\Delta_0$ coincides with the canonical coproduct of the quantum subgroup $\clu/\cli_0$ inherited from that of $\clu$. It is also easy to see that the object $(\widetilde{\clg}, \Delta_0, \Gamma_{\cli_0})$ is universal in the category ${\bf Q}^\prime_R$, using the fact that (by Lemma \ref{lem2}) any compact quantum group $(\cls, U)$ acting by volume and orientation preserving isometries on the given spectral triple is isomorphic with a quantum subgroup $\clu/\cli$, for some Hopf $C^*$-ideal $\cli$ of $\clu$. 

Finally, the faithfulness of $U_0$ follows from the universality by standard arguments which we briefly sketch. If $\clg_1 \subset \widetilde{\clg}$ is a $\ast$-subalgebra of $\widetilde{\clg}$ such that $\widetilde{U_0} \subseteq \clm(\clk(\clh) \ot \clg_1)$, it is easy to see that $(\clg_1, U_0)$ is also a universal object, and by definition of universality of $\widetilde{\clg}$ it follows that there is a unique morphism, say $j$, from $\widetilde{\clg}$ to $\clg_1$. But the map $i \circ j$ is a morphism from $\widetilde{\clg}$ to itself, where $i : \clg_1 \raro \widetilde{\clg}$ is the inclusion. Again by universality, we have that $i \circ j ={\rm id}_{\widetilde{\clg}}$, so in particular, $i$ is onto, i.e. $\clg_1=\widetilde{\clg}$.  
\qed

Consider the  $\ast$-homomorphism $\alpha_0:=\alpha_{U_0}$, where  $(\widetilde{\clg},U_0)$ is the universal object obtained in the previous theorem. For every state $\phi$ on $\widetilde{\clg}$, $({\rm id} \ot \phi) \circ \alpha_0$ maps $\cla$ into $\cla^{\prime \prime}$. However, in general $\alpha_0$ may not be faithful  even if $U_0$ is so, and let $\clg$ denote the $C^*$-subalgebra of $\widetilde{\clg}$ generated by the elements $\{ (f \ot {\rm id})\circ \alpha_0(a), f \in \cla^*, a \in \cla \}$. 

\brmrk

If the spectral triple is even, then all the proofs above go through with obvious modifications.

\ermrk

\bdfn
We shall call $\clg$ the quantum  group of orientation -preserving isometries of the R- twisted spectral triple $(\cla^\infty, \clh, D, R)$  and denote it by $ {QISO^+}_R ( \cla^{\infty}, \clh, D, R ) $ or even simply as ${QISO^{+}}_{R}(D)$. The quantum group $\widetilde{\clg}$ is denoted by $\widetilde{QISO^{+}}_{R}(D)$. 

If the spectral triple is even, then we will denote 
$ \clg $ and $ \widetilde{\clg} $ by $ QISO^{+}_{R} ( D, \gamma ) $ and  $ \widetilde{QISO^{+}_{R}} ( D, \gamma ) $ respectively.

\edfn

\subsection{Stability and topological action}
it is not clear from the definition and construction of $QISO_R^+(D)$ whether the $C^*$ algebra $\cla$ generated by $\cla^\infty$ is stable under $\alpha_0$   in the sense that $({\rm id} \ot \phi) \circ \alpha_0$ maps $\cla$ into $\cla$ for every $\phi$. Moreover, even if $\cla$ is stable, the question remains whether $\alpha_0$ is a $C^*$-action of the CQG $QISO^+_R(D)$. Although we could not yet decide whether the general answer to the above two questions are afirmative, we should point out that for all the explicit examples considered by us so, both these questions do have an affirmative answer. It is perhaps tempting (though a bit too daring at the moment) to conjecture that in general, $QISO^+_R$ will have a topological, i.e. $C^*$ action. Indeed, we have quite a few results in this direction. Let us list the cases where we are able to prove that $QISO^+_R(D)$ has $C^*$-action:\\
(i) For any spectral triple for which there is a 'reasonable' Laplacian in the sense of \cite{goswami}. This includes all classical spectral triples as well as their Rieffel deformation ( with $R=I$ ).\\
(ii) Under the assumption that there is an eigenvalue of $D$ with a one-dimension eigenspace spanned by a cyclic separating vector $\xi$ such that any 
eigenvector of $D$ belongs to the span of $\cla^\infty \xi$ and $\{ a \in \cla^\infty:~a \xi ~{\rm is}~{\rm an}~{\rm eigenvector}~{\rm of}~D \}$ is norm-dense in $\cla$ ( to be proved in subsection 2.4).\\
(iii) Under some analogue of the classical Sobolev conditions w.r.t. a suitable group action on $\cla$ ( see \cite{goswami_action}). \\

.
Now we prove the sufficient conditions in (i).

We begin with a sufficient condition for stability of $\cla^\infty$ under $\alpha_0$.
 Let $ ( \cla^{\infty}, \clh, D ) $ be a (compact type) spectral triple such that \vspace{1mm}\\
({\bf 1}) ${\cla^\infty}$ and $\{ [D,a], ~a \in {\cla^\infty} \}$ are 
contained in the domains of all powers of the derivations $[ D, \cdot ] $ and $ [|D|, \cdot] .$ \\

We will denote by $ \widetilde{T_t} ,$ the one parameter group of $ \ast $ automorphisms on $ \clb ( \clh ) $ given by $ \widetilde{T_t} ( S ) =
 e^{ it D} S e^{-it D} ~ \forall S \in \clb ( \clh ).$ We will denote the generator of this group by $ \delta.$  
For $ X $ such that $ [ D, X ] $ is bounded, we have $ \delta ( X ) = i[ D, X ]  $ and hence
$${\left\| \widetilde{T_t} ( X ) - X \right\|} = {\left\| \int^{t}_{0} \widetilde{T_s} ([ D, X ])ds \right\|} \leq ~ t {\left\| [ D, X ] \right\|} .$$

Let us say that the spectral triple satisfies the {\it Sobolev condition} if $$\cla^\infty=\cla^{\prime \prime} \bigcap_{n \geq 1} {\rm Dom}(\delta^n).$$  Then we have the following result, which is a natural generalization of the classical situation, where a measurable isometric action automatically becomes topological (in fact smooth).
\bthm
\label{stable}
(i) For every state $\phi$ on $\clg$, we have  $({\rm id} \ot \phi)\circ \alpha_0(\cla^\infty) \in \cla^{\prime \prime} \bigcap_{n \geq 1} {\rm Dom}(\delta^{n})$. 
 (ii) If the spectral triple satisfies the  Sobolev condition then  
  $\cla^\infty$ ( and hence $\cla$ ) is stable under $\alpha_0$.   
\ethm
{\it Proof:}\\
Since $U_0$ commutes with $D \ot I$, it is clear that the automorphism group $\widetilde{T_t}$ commutes with $\alpha_0^\phi \equiv ({\rm id} \ot \phi) \circ \alpha_0$, and thus by the continuity of $\alpha_0$ in the stromg operator topology it is easy to see that, for $a \in {\rm Dom}(\delta),$ we have \bean \lefteqn {\lim_{t \raro 0+} \frac{\widetilde{T_t}(\alpha_0^\phi(a))-\alpha_0^\phi(a)}{t}}\\
&=&   \lim_{t \raro 0+} \alpha_0^\phi(\frac{\widetilde{T_t}(a)-a}{t})\\
&=& \alpha_0^\phi(\delta(a)).\eean
Thus, $\alpha_0^\phi$ leaves ${\rm Dom}(\delta)$ invariant and commutes with $\delta$. Proceeding similarly, we prove (i). The assertion (ii)  is a trivial consequence of (i) and the  Sobolev condition. \qed

Let us now assume that \\
({\bf 2}) The spectral triple is $\Theta$-summable, i.e. for every $t>0$, $e^{-tD^2}$ is trace-class.\\
  Consider the functional  $ \tau ( X ) = 
{\rm Lim}_{t \rightarrow 0} \frac{ {\rm Tr} ( X e^{- t D^2} )}{ {\rm Tr} ( e^{- t D^2} )}$ is a (not necessarily faithful or normal)
 state on $ \clb ( \clh ) $ where ${\rm Lim}$ is a suitable Banach limit as in \cite{fro}.  
Moreover, $ \tau $ is a positive faithful trace on the $ \ast $ algebra, say $\cls^\infty$, generated by $ \cla^{\infty} $ and $ \{ [ D, a ] : a \in \cla^{\infty}\}, $
 which is to be interpreted as the volume form( we refer to \cite{fro}, \cite{goswami} for the details ). The completion  of $\cls^\infty$ in the norm of $\clb(\clh)$ is denoted by $\cls$.

From the definition of $ \tau ,$ it is also clear that $ \widetilde{T_t} $ preserves $ \tau ,$ so extends to a group of unitaries on $\cln:=L^2( \cls^\infty, \tau)$.
 Moreover,  for $X$ such that $[D,X] \in \clb(\clh)$, in particular for $X \in \cls^\infty$, we  have (denoting by $\| a \|_2$ the $L^2$-norm $\tau(a^*a)^{\frac{1}{2}}$  and the operator norm of $\clb(\clh)$ by $\| \cdot \|_\infty$) 
\bean \lefteqn{\left\| \widetilde{T_s} ( X ) - X  \right\|^2_2 }\\
& =& \tau ( X^* ( X - \widetilde{T_s} ( X ) ) ) + \tau ( {( X - \widetilde{T_s} ( X ) )}^{*} ) X ) \\
  &  \leq &   2  \left\| X - \widetilde{T_s} ( X )  \right\|_{\infty} \left\| X \right\|_2 \\
& \leq & s \| [D,X] \|_\infty \| X \|_2,\\
\eean
which clearly shows that $s \mapsto \widetilde{T_s}(X)$ is $L^2$-continuous for $X \in \cls^\infty$, hence (by unitarity of $\widetilde{T_s}$) on the whole of $\cln$, i.e. it is a strongly continuous one-parameter group of unitaries. Let us denote its generator by $\tilde{\delta}$, which is a skew adjoint map, i.e. $i \tilde{\delta}$ is self adjoint, and $\widetilde{T_t}={\rm exp}(t \tilde{\delta})$. Clearly, $\tilde{\delta}=\delta=[D, \cdot]$ on $\cls^\infty$.   

We will denote $ L^{2}  ( \cla^\infty, \tau ) \subset \cln$ by $ \clh^{0}_{D}$ and the restriction of $\tilde{\delta}$ to $\clh^0_D$ (which is a closable map from $\clh^0_D$ to $\cln$) by $d_D$. Thus, $d_D$ is closable too. 

  .

We now recall the assumptions made in \cite{goswami}, for definining  the `Laplacian' and the corresponding  quantum isometry group of a spectral triple $(\cla^\infty,\clh,D)$, without going into all the technical details, for which the reader is referred to \cite{goswami}. 

The following conditions will also be assumed throughout the rest of this subsection:\vspace{1mm}\\

({\bf 3}) $ \cla^{\infty} \subseteq {\rm Dom} ( \cll ) $ where $ \cll \equiv \cll_D := - d^{*}_{D} d_{D} .$\\
({\bf 4}) $ \cll $ has compact resolvent.\\
({\bf 5}) Each eigenvector of $ \cll $ ( which has a discrete spectrum , hence a complete set of eigenvectors ) belongs to $ \cla^{\infty}.$\\
({\bf 6}) The complex linear span of the eigenvectors of $ \cll $, say $ \cla^{\infty}_{0} $ ( which is a subspace of $ \cla^{\infty} $ by assumption ({\bf 5.})  ), is norm dense in $ \cla^{\infty}.$\vspace{2mm}\\

It is clear that  $\cll(\cla^\infty_0) \subseteq \cla^\infty_0$. The $\ast$-subalgebra  of $\cla^\infty$ generated by $\cla^\infty_0$ is denoted by $\cla_0$.
  We also note that $\cll=P_0\tilde{\cll} P_0$, where $\tilde{\cll}:= (i \tilde{\delta})^2$ (which is a self adjoint operator on $\cln$) and $P_0$ denotes the orthogonal  projection in $\cln$ whose range is  the subspace $\clh^0_D$.

\bthm

\label{QISO_I_D_<_Q_L_D}
 
 Let $( \cla^{\infty}, \clh, D )$ be a spectral triple satisfying the assumptions $({\bf 1})-({\bf 6})$ made above. In addition, assume that  at least one of conditions (a) and (b) mentioned below is satisfied:\\ 
 
(a) $ {\cla}^{\prime \prime} \subseteq \clh^{0}_{D} .$\\
(b) $ \alpha_0^{\phi} ( \cla^{\infty} ) \subseteq \cla^{\infty} $ for every state $\phi$ on $\clg=QISO^+_R(D)$.\\
 
 Then  $ \alpha_0 $ is a $C^*$-action of $ QISO^{+}_{I} ( D ) $  on $ \cla .$ 

\ethm

{\it Proof :}\\
Under either of the conditions (a) and (b), the map $\alpha^\phi_0$ maps $\cla^\infty$ into $\clh^0_D \subseteq \cln$ (for any fixed $\phi$). Since $\alpha^\phi_0$ also commutes with $[D, \cdot]$ on $\cla^\infty$, it is clear that $\alpha^\phi_0(\cls^\infty) \subseteq \cln$ too. In fact, 
 using the complete positivity of the map $ \alpha_0^{\phi} $ and the $\alpha_0$-invariance of $\tau$, we see that
 $$ {\tau ( \alpha_0^{\phi} ( a )}^{*} \alpha_0^{\phi} ( a ) ) \leq \tau ( \alpha_0^{\phi} ( a^* a ) )
 = ( {\rm id} \otimes \phi ) ( \tau \otimes {\rm id} ) \alpha _0( a^* a ) = \tau( a^* a ).1, $$ which implies that $ \alpha_0^{\phi} $ extends to a bounded operator from $ \cln$ to itself.  Since $U_0$ commutes with $D$, it is clear that $\alpha_0^\phi$ (viewed as a bounded operator on $\cln$) 
 will commute with the group of unitaries $\widetilde{T_t}$, hence with its generator $\tilde{\delta}$ and also with the self adjoint operator $\tilde{\cll}=(i \tilde{\delta})^2$. 

On the other hand, it follows from the  definition of $  \clg=QISO^{+}_{I}( D )$ that $(\tau \ot {\rm id})(\alpha_0(X))=\tau(X) 1_\clg $ for all 
$X \in \clb(\clh)$, in particular for $X \in \cls^\infty$, and thus the map $\cls^\infty \ot_{\rm alg} \clg \ni (a \ot q) \mapsto \alpha_0(a)(1 \ot q)$ 
extends to a $\clg$-linear unitary, denoted by $W$ (say), on the Hilbert $\clg$-module $\cln \ot \clg$. 
Note that here we have used the fact that for any $\phi$, $({\rm id} \ot \phi)(W) (\cls^\infty \ot_{\rm alg} \clg) \subseteq \cln$, 
since $\alpha^\phi_0(\cls^\infty) \subseteq \cln$. The commutativity of $\alpha^\phi_0$ with $\widetilde{T_t}$ for every $\phi$ clearly implies that 
$W$ and $\widetilde{T_t} \ot {\rm id}_\clg$ commute on $\cln \ot \clg$. Moreover, $\alpha^\phi_0$ maps $\clh^0_D$ into itself, so $W$ maps 
$\clh^0_D \ot \clg$ into itself, and  hence (by unitarity of $W$) it commutes with the projection $P_0 \ot 1$. 
It follows that $\alpha^\phi_0$ commutes with $P_0$, and (since it also commutes with $\tilde{\cll}$), hence commutes with $\cll=P_0 \tilde{\cll} P_0$ 
as well.

Thus, $ \alpha^\phi_0 $ preserves each of  the (finite dimensional)  eigenspaces of the Laplacian $\cll$,  and so is a Hopf algebraic action on  the subalgebra $\cla_0$ spanned algebraically by these eigenspaces. 
Moreover, the   $\clg$-linear unitary $W$ clearly restricts to a unitary representation on each of the above eigenspaces. 
If we denote by  $ (( q_{ij} ))_{( i,j )} $ the $\clg$-valued unitary matrix corresponding to one  such particular eigenspace, then by the general theory of CQG representations, $q_{ij}$ must belong to $\clg_0$ and we must have $ \epsilon ( q_{ij} ) = \delta_{ij} $ (Kronecker delta).  This  implies $ ( {\rm id} \otimes \epsilon ) \circ \alpha_0  = {\rm id} $ on each of the eigenspaces, hence on the norm-dense subalgebra $\cla_0$ of $\cla$,  completing the proof of the fact that $ \alpha_0 $ extends to a $ C^* $ action on $ \cla. $ \qed

 Combining the above theorem with Theorem \ref{stable}, we get the following immediate corollary.
\bcrlre
If the spectral triple satisfies the Sobolev condition mentioned before, in addition to the assumptions ${\bf 1}-{\bf 6}$, then $QISO^+_I(D)$ has a $C^*$-action.  In particular, for a classical spectral triple, $QISO^+_I(D)$ has $C^*$-action.
\ecrlre
\brmrk 
Let us remark here that in case the restriction of $\tau$ on $\cla^\infty$ is normal, i.e. continuous w.r.t. the WOT inherited from $\clb(\clh)$, 
then $\clh^0_D$ will contain $\cla^{\prime \prime}$, which is the WOT closure of $\cla^\infty$ in $\clb(\clh)$, .e. condition (a) of Theorem \ref{QISO_I_D_<_Q_L_D} (and hence its conclusion) holds.  
\ermrk

\brmrk
The results obtained in this subsection can be formulated and proved in an $R$-twsited set-up as well, if the corresponding Laplacian (which is an extension of $d_{D,R}^* d_{D,R}$, where $d_{D,R}^*$ denotes the adjoint of $d_D \equiv d_{D,R}$ w.r.t. the $R$-twisted volume form) `exists' and satisfies the analogues of the assumptions made in this subsection about $\cll_D$. In \cite{sphere}, we have made some computations with such an $R$-twisted Laplacian arising naturally in that context.
\ermrk

\brmrk
In a private communication to us, Shuzhou Wang has kindly pointed out  that a possible alternative approach  to the formulation of quantum group of isometries may involve the category of CQG which has a $C^*$-action on the underlying $C^*$ algebra and a unitary representation w.r.t. which the Dirac operator is equivariant. We are not sure whether or how  it is  possible to show the existence of a universal object (even after fixing a choice of volume form) in this category; however, if the existence can be proved then the universal object will automatically have $C^*$ action.
\ermrk

\subsection{ Universal object in the  categories  ${\bf Q} $ or ${\bf Q}^{\prime}$}

We shall now investigate   further conditions  on the spectral triple which will ensure the existence of a  universal object in the category ${\bf Q}$ or  ${\bf Q}^\prime$.
 Whenever scuh an universal object exists we shall denote it by $\widetilde{QISO^+}(D)$, and denote by $QISO^+(D)$ its  largest Woronowicz subalgebra  for which $\alpha_U$ on $\cla^\infty$  (where $U$ is the unitary representation of $\widetilde{QISO^+}(D)$ on $\clh$) is faithful.
 \brmrk
If ${\widetilde{QISO^+}}(D)$ exists, by \cite{goswami_rmp}, there will exist some $R$ such that ${\widetilde{QISO^+}}(D)$ is an object in the category ${\bf Q}^\prime_R(D)$. Since the universal object in this category, i.e. $\widetilde{QISO^+}_R(D)$, is clearly a sub-object of $\widetilde{QISO^+}(D)$, we have $\widetilde{QISO^+}(D) \cong \widetilde{QISO^+}_R(D)$ for this choice of $R$.
\ermrk

  Let us state and prove a result below, which gives some sufficient conditions for the existence of $\widetilde{QISO^+}(D)$.
\bthm

\label{unrestricted}

Let $(\cla^\infty, \clh, D)$ be a spectral triple of compact type as before and assume that $D$ has an one-dimensional eigenspace spanned by a unit vector $\xi$, which is cyclic and separating for the algebra $\cla^\infty$. Moreover, assume that each eigenvector of $D$ belongs to the dense subspace $\cla^\infty \xi$ of $\clh$. Then there is a universal object, 
$(\widetilde{\clg}, U_0)$. Moreover, $\widetilde{\clg}$ has a coproduct $\Delta_0$ such that $(\widetilde{\clg},\Delta_0)$ is a compact quantum group and $(\widetilde{\clg},\Delta_0,U_0)$ is a universal object in the category ${\bf Q}^\prime.$ 
 
 We will denote by $ \clg $ the Woronowicz $ C^{\ast} $ subalgebra of $ \widetilde{\clg} $ generated by elements of the form $ \left\langle ( \xi \otimes 1 ) \alpha_{U_{0}} ( a ) , ~ \eta \otimes 1 \right\rangle $ where $ \xi, \eta \in \clh, a \in \cla^{\infty} $ and $ \left\langle . ,~ . \right\rangle $ denotes the $ \widetilde{\clg} $ valued inner product of $ \clh \otimes \widetilde{\clg} .$ We have $ \widetilde{\clg} \cong \clg \ast C ( \IT ).$
\ethm
{\it Proof:}\\
Let $V_i$, $\{ e_{ij} \}$ be as before, and by assumption $e_{ij}=x_{ij} \xi$ for a unique $x_{ij} \in \cla^\infty$. Clearly, since $\xi$ is separating, the vectors $\{ \overline{e_{ij}}=x^*_{ij} \xi, j=1,..., d_i \}$ are linearly independent, so the matrix $Q_i=(( \lgl \overline{e_{ij}}, \overline{e_{ik}} ))_{j,k=1}^{d_i}$ is positive and invertible. Now, given a quantum family of orientation-preserving  isometries $(\cls, U)$, we must have $\widetilde{U} ( \xi \ot 1)= \xi \ot q$, say, for some $q$ in $\cls$, and from the unitarity of $\widetilde{U}$ it follows that $q$ is a unitary element. Moreover, $U$ leaves $V_i$ invariant, so let $\widetilde{U}(e_{ij} \ot 1)=\sum_k e_{ik} \ot v^{(i)}_{kj}$. But this can be rewritten as $$ \alpha_U(x_{ij}) (\xi \ot q)=\sum_k x_{ik} \xi \ot v^{(i)}_{kj}.$$ Since $\xi$ is separating and $q$ is unitary, this implies $ \alpha_U(x_{ij})=\sum_k x_{ik} \ot v^{(i)}_{kj}q^*,$  and thus we have $$ \widetilde{U}(\overline{e_{ij}}\ot 1)=\alpha_U(x_{ij})^*(\xi \ot q)=\sum_k x^*_{ik} \xi \ot q^2 (v^{(i)}_{kj})^*=\sum_k \overline{e_{kj}} \ot q^2 (v^{(i)}_{kj})^*.$$ 
  Taking the $\cls$-valued inner product $\lgl \cdot, \cdot \rgl_\cls$ on both sides of the above expression,   and using the fact that $U$ preserves this $\cls$-valued inner product,  we obtain $Q_i=v_i^\prime Q_i \overline{v_i}$ (where $v_i=(( v^{(i)}_{kj}))$). Thus, $Q_i^{-1}v_i^\prime Q_i$ must be the (both-sided) inverse of $\overline{v_i}$. Thus, we get a canonical surjective morphism from $ A_u(Q_i) $ to the $C^{\ast}$ algebra generated by $v^{(i)}_{kj}$. This induces a surjective morphism from the free product of $A_u(Q_i)$, $i=1,2,...$ onto $\cls$. The rest of the arguements for showing the existence of $\widetilde{\clg}$  will be quite similar to the arguements used in the proof of Theorem \ref{main}, hence omitted. It is also quite obvious from the proof that $\widetilde{\clg} =\clg \ast C^*(q) \cong \clg \ast C(\IT)$.
\qed

\brmrk

Some of the examples considered in the next section will show that the conditions of the above theorem are not actually necessary; $\widetilde{QISO}^+(D)$ may exist without the existence of a single cyclic seperating eigenvector as above.  
\ermrk

Let $ ( {\cla}^{\infty}, \clh, D ) $ be a spectral triple satisfying the conditions of the above theorem.

Let the faithful vector state corresponding to the cyclic separating vector $ \xi $ be denoted by $ \tau.$
Let $ \cla^{\infty}_{0} = {\rm span} \{ a \in \cla^{\infty} : a \xi $ is an eigenvector of $ D  \} $

Moreover, assume that $ \cla^{\infty}_{0} $ is norm dense in $ \cla^{\infty} .$

Let $ \hat{D} : \cla^{\infty}_{0} \rightarrow  \cla^{\infty}_{0} $ be defined by :

$ \hat{D} ( a ) \xi = D ( a \xi ) $

This is well defined since $ \xi $ is cyclic separating.

\bdfn
Let $ \cla $ be a $ C^{\ast} $ algebra and $ {\cla}^{\infty} $ be a dense $ \ast $ subalgebra.
Let $ ( {\cla}^{\infty} , \clh, D ) $ be a spectral triple as above.

Let $ \widehat{\bf C} $ be the category with objects $ ( \clq, \alpha ) $ such that $ \clq $ is a compact quantum group with an action $\alpha$  on $ \cla $ such that:

1. $ \alpha $ is $ \tau $ preserving ( where $ \tau $ is above ), i.e, $ ( \tau \otimes id ) \alpha ( a ) = \tau ( a ).1 $

2. $ \alpha $ maps $ \cla^{\infty}_{0} $ inside $  \cla^{\infty}_{0} \otimes_{alg} \clq .$

3. $ \alpha \widehat{D} = ( \widehat{D} \otimes I ) \alpha  .$

\edfn

\bcrlre

\label{unrestrictedcorollary}
 
 There exists a universal object $ \widehat{\clq} $ of the category $ \widehat{\clc} $ and it is isomorphic to the Woronowicz $ {C}^{*} $ subalgebra $\clg={QISO}^+(D)$ of $ \widetilde{\clg} $  obtained in Theorem \ref{unrestricted}.

\ecrlre

{\it Proof :}\\ The proof of the  existence of the universal object follows verbatim from the proof of Theorem 4.7 in \cite{goswami}  replacing $ \cll $ by $ \widehat{D} $ and noting that $ D $ has compact resolvent. We denote by $\widehat{\alpha} $ the action of $ \widehat{\clq} $ on $ \cla.$ 

Now, we prove that $ \widehat{\clq} $ is isomorphic to  $ \clg .$

 Each eigenvector of $ D$ is in  $  \cla^{\infty} $ by assumption. It is easily observed from the proof of Theorem \ref{unrestricted} that  $ \alpha_{U_{0}}$ maps the norm-dense $\ast$-subalgebra $  {\cla_{0}}^{\infty} $ into $ {\cla_{0}}^{\infty} \otimes_{alg} \clg_{0} $, and $({\rm id} \ot \epsilon)\circ \alpha_{U_0}={\rm id}$, so that $\alpha_{U_0}$ is indeed an action of the CQG $\clg$.   Moreover, it can be easily seen that $ \tau $ preserves  $  \alpha_{U_{0}} $ and that $ \alpha_{U_{0}} $ commutes with $ \widehat{D}.$ Therefore,$ ( \clg, \alpha_{U_{0}} ) \in {\rm Obj} ( \widehat{\bf C} ) $,  and  hence  $ \clg $ is a quantum subgroup of $ \widehat{\clq} $ by the universality of $\widehat{\clq}$.

For the converse, we start by showing that $ \widehat{\alpha} $ induces a unitary representation $W$ of $ \widehat{\clq} \ast C( \IT ) $ on $ \clh $ which commutes with $ D$, and the corresponding conjugated action $\alpha_W$ coincides with $\widehat{\alpha}$. 

Define $ W( a \xi ) = \widehat{\alpha}( a ) ( \xi ) ( 1 \otimes q^{*} )  \forall a \in \cla^{\infty}_{0} $ where $ q  $ is a generator of $ C( \IT ).$

Since we have $ ( \tau \otimes id ) \alpha ( a ) = \tau ( a ).1,$ it follows that $ \widetilde{W}  $ is a ($\widehat{\clq} \ast C(\IT)$-linear)  isometry on the dense subspace $ \cla^{0}_{\infty}\xi \ot_{\rm alg} \widehat{\clq} $ and thus extends to $ \clh \ot \widehat{\clq}\ast C(\IT)$ as an isometry. Moreover, since $ \widehat{\alpha}( \cla ) ( 1 \otimes \widehat{\clq} )$ is norm dense in $ \cla \otimes \widehat{\clq}$ (by the definition of a CQG action)  it is clear that the range of $\widetilde{W}$ is dense, so $\widetilde{W}$ is indeed a unitary.  It is quite obvious that it is a unitary representation of $\widehat{\clq} \ast C(\IT)$.

We also have, \bean 
\lefteqn{W D ( a \xi ) }\\
& =&  W( \widehat{D} ( a ) \xi ) =
   \widehat{\alpha} ( \widehat{D} ( a ) ) ( \xi ) ( 1 \otimes q^{*} ) \\
& =&  ( D \otimes I ) ( \widehat{\alpha} ( a ) \xi ) ( 1 \otimes q^{*} ) 
 = ( D \otimes I ) W ( a \xi ) , \eean

i.e.   $ W $ commutes with $ D .$

It is also easy to observe that $ \alpha_{W} = \widehat{\alpha}.$ This gives a surjective CQG morphism from $ \widetilde{\clg} = \clg \ast C( \IT ) $ to $ \widehat{\clq} \ast C( \IT ),$ sending $ \clg $ onto $ \widehat{\clq},$ which completes the proof.

\qed

\section{Comparison with the approach of \cite{goswami} based on Laplacian}
Throughout this section, we shall assume the set-up of  subsection 2.3 for the existence of a `Laplacian',  including assumptions ${\bf 1}-{\bf 6}$. Let us also use the notation of that subsection. 

 As in \cite{goswami}, we say that a CQG $(\cls,\Delta)$ which has an action $\alpha$ on $\cla$ is said to act smoothly and isometrically on the noncommutative manifold $(\cla^\infty, \clh, D)$ if $({\rm id} \ot \phi) \circ \alpha(\cla^\infty_0) \subseteq \cla^\infty_0$ for all state $\phi$ on $\cls$, and also $({\rm id} \ot \phi) \circ \alpha$ commutes with the Laplacian $\cll \equiv \cll_D$ on $\cla^\infty_{0}$. One can consider the category $ {\bf Q}^{\cll_D} $ of all compact quantum groups acting smoothly and isometrically on $\cla$, where the morphisms are quantum group morphisms which intertwin the actions on $\cla$. 
We make the following additional assumption throughout the present section:\vspace{2mm}\\
 ({\bf 7}) There exists a universal object in $ {\bf Q}^{\cll_D} $ (i.e. the quantum isometry group for the Laplacian $\cll \equiv \cll_D$ in the sense of \cite{goswami}), and it is denoted by  $ \clq^\cll \equiv \clq^{\cll_{D}} $ \vspace{2mm}\\

\brmrk
It is proved in \cite{goswami}  that that under the additional `connectedness assumption ' ( the kernel of $ \cll $ is one dimensional, spanned by the identity $1$ of $ \cla^{\infty} $ viewed as a unit vector in $ \clh^{0}_{D} $), the category ${\bf Q}^\cll$ has a universal object, say $\clq^\cll$, called the quantum isometry group in \cite{goswami}.  
In \cite{goswami} it was also shown ( Lemma 2.5, $ b \Rightarrow a $ ) that for an isometric group action on a  classical manifold ( not necessarily connected ), the  volume functional is automatically preserved. It can be easily seen that the proof goes verbatim for a quantum group action. As this volume preserving condition implies the existence of $ \clq^{\cll_{D}},$ for a classical manifold ( not necessarily connected ) $ \clq^{\cll_D} $ always exists. 

\ermrk

The following result now follows immediately from Theorem \ref{QISO_I_D_<_Q_L_D} of subsection 2.3. 

\bcrlre

\label{QISO_I_D_<_Q_L}

$ QISO^+_{I}( D ) $ is a sub-object of $ \clq^{\cll_{D}}  $ in the category $ {\bf Q}^{\cll_D} .$

\ecrlre
{\it Proof:}\\
The proof is a consequence of the fact that $QISO^+_I(D)$ has the $C^*$-action $\alpha_0$ on $\cla$, and the observation already made in the proof of  the Theorem \ref{QISO_I_D_<_Q_L_D} that this action commutes the Laplacian $\cll_D$. 
\qed 

We recall the Hilbert space of forms from \cite{fro}.
Given the spectral triple $ ( \cla^{\infty}, \clh, D )$, consider the Hilbert space $ \clh_{d + d^*} = \oplus \clh^{n}_{D} ,$
( we refer to \cite{fro} for the definition and other details  ) and the densely defined self adjoint operator $ d + d^* $ on it. 
There is a representation of $ \cla $ on this Hilbert space given by $ \pi_{d + d^*} ( a ) ( [ a_0 d_{D} a_1.....d_{D} ( a_k ) ] ) =
 [ a a_0 d_{D} a_1.....d_{D} ( a_k ) ].$ where $ [  a_0 d_{D} a_1.....d_{D} ( a_k ) ] $ denote the class of  $ a_0 d_{D} a_1.....d_{D} ( a_k ) $ 
in $ \Omega^{k}_{D} ( {\cla}^{\infty} ) = \pi_{D} ( \Omega_{k} ( \cla^{\infty} ) ) / J_k + \delta ( J_{k - 1} ) $ 
where $ \delta, ~ J_k $ is as in Page 124, \cite{fro}.  Then $( \cla^{\infty}, \clh_{d + d^*}, d + d^* ) $ is another spectral triple. We assume that this 
 is of compact type, i.e. $d+d^*$ has compact resolvents.

We will denote the inner product on the space of $ k $ forms coming from the spectral triples $ ( \cla^{\infty}, ~ \clh, ~ D ) $ and 
$( \cla^{\infty}, ~ \clh_{d + d^*},~ d + d^* ) $ by $ {\left\langle ~ , ~ \right\rangle}_{\clh^{k}_{D}} $ and 
$ {\left\langle ~ , ~ \right\rangle}_{\clh^{k}_{d + d^*}} $ respectively, $ k = 0,1.$ 

We will denote by $ \pi_{D} , ~ \pi_{d + d^*} $ the representations of $ \cla^{\infty} $ in $ \clh $ and $ \clh_{d + d^*} $ respectively.

Let $ U_{d + d^*} $ be the canonical unitary representation of 
$ QISO^{+}_{I}( d + d^* )$ on $ \clh_{d + d^*}.$ 

$ \clh_{d + d^*} $ breaks up into finite dimensional orthogonal subspaces corresponding to the distinct eigenvalues of $ \Delta:=(d + d^*)^2=d^*d+dd^*.$ 
It is easy to see that $\Delta$ leaves each of the subspaces $\clh^i_D$ invariant, and we will denote by $ V_{\lambda,i} $ the subspace of $ \clh^{i}_{d+d^*} $ spanned by eigenvectors of $ \Delta $ corresponding to the eigenvalue 
$ \lambda.$ Let $ \{ e_{j,\lambda,i} \}_{j} $ be an orthonormal basis of $ V_{\lambda,i}.$ Note that $\cll_D$ is the restriction of $\Delta$ to $\clh^0_D$.

Now we recall the result of Section 2.4 of \cite{goswami}. It was shown there that $ \clq^{\cll_{D}} $ has a unitary representation $ U \equiv U_\cll $ on
 $ \clh^{d + d^*} $ such that $ U $ commutes with $ d + d^* .$ Thus $ ( \cla^{\infty}, ~  \clh_{d + d^*}, ~ d + d^* ) $ is a $ \clq^{\cll_{D}} $ 
equivariant spectral triple. It follows from the construction in \cite{goswami} that $ \clq^{\cll_{D}} $ is a quotient ( by a Woronowicz $ C^* $ ideal ) of the free product of countably many Wang algebras of the type  $ A_{u} ( I )  $,  and hence is a compact  quantum Kac algebra. Thus, it has tracial Haar state, which implies, by Theorem 3.2 of \cite{goswami_rmp}, that  $ \alpha_U $ keeps the functional $ \tau_I $ invariant. Thus, we have:

\bppsn

\label{Q_L_leq_QISO(d + d*)}

 $ ( \clq^{\cll_{D}}, ~ U_\cll ) $ is a sub object of $ ( QISO^{+}_{I} ( d + d^* ), ~ U_{d + d^*} ) $ in the category $ {\bf Q}_{I} ( d + d^* ),$ so in particular,
  $ \clq^{\cll_{D}} $ is  isomorphic to a quotient of $ QISO^{+}_{I}( d + d^* ) $ by a Woronowicz $ C^* $ ideal.

 \eppsn

We shall  give (under mild conditions) a concrete description of the above Woronowicz ideal. 

Let $ \cli $ be the $ C^* $ ideal of $QISO^{+}_{I}( d + d^* )$ generated by

 $$ \cup_{\lambda \in \sigma( \Delta )} \{ \left\langle ( P^{\bot}_{0} \otimes id ) U_{d + d^*} ( e_{j \lambda 0} ), e_{j \lambda i^{\prime}} 
\otimes 1 \right\rangle : j, i^{\prime} \geq 1 \},$$
where $ P_0 $ is the projection onto $ \clh^{0}_D $ and $ \left\langle . ,~ . \right\rangle $ denotes the $ QISO^{+}_{I}( d + d^* ) $ 
valued inner product.

Since $ U_{d + d^*} $ keeps the eigenspaces of $ \Delta=(d + d^*)^2 $ invariant, we can write 
$$ U ( e_{j \lambda 0} ) = \sum_{k} e_{k \lambda 0} \otimes ~ q_{k j \lambda 0} + \sum_{i^{\prime}\neq 0, k^{\prime} } e_{k^{\prime} \lambda i^{\prime}} 
\otimes ~ q_{k^{\prime} j \lambda i^{\prime}}, $$ for some $ q_{k j \lambda 0}, ~ q_{k^{\prime} j \lambda i^{\prime}}  \in QISO^{+}_{I}( d + d^* ).$

We note that $ q_{k^{\prime} j \lambda i^{\prime}} \in \cli $ if $ i^{\prime} \neq 0.$

\blmma

$ \cli $ is a co-ideal of $ QISO^{+}_{I} ( d + d^* ).$

\elmma

{\it Proof :}\\ It is enough to prove the relation $ \Delta ( \cli ) \in \cli \otimes QISO^{+}_{I}( d + d^* ) + QISO^{+}_{I}( d + d^* ) \otimes \cli $ for the elements in $ \cli $ of the form $ \left\langle ( P^{\bot}_{0} \otimes id ) U_{d + d^*} ( e_{j \lambda  0} ), e_{j \lambda i_0} \otimes 1 \right\rangle .$
We have:
\bean \Delta ( \left\langle ( P^{\bot}_{0} \otimes id ) U_{d + d^*} ( e_{m \lambda 0} ), e_{j \lambda i_0} \otimes 1 \right\rangle ) \eean

 \bean = \left\langle ( P^{\bot}_{0} \otimes id ) ( id \otimes \Delta ) U_{d + d^*} ( e_{m \lambda 0} ), e_{j \lambda i_0} \otimes 1 \ot 1 \right\rangle ) \eean
 
 \bean = \left\langle ( P^{\bot}_{0} \otimes id ) U_{(12)} U_{(13)} ( e_{m \lambda i} ), e_{j \lambda i_0} \otimes 1 \ot 1 \right\rangle  \eean
 
 \bean = \left\langle ( P^{\bot}_{0} \otimes id ) U_{(12)} ( \sum_{k} e_{k \lambda 0} \otimes 1 \otimes q_{k m \lambda 0} ) ~~, ~~ e_{j \lambda i_{0}} \otimes 1 \ot 1 \right\rangle  \eean 
 
\bean +  \sum_{i^{\prime} \neq 0,l} \left\langle ( P^{\bot}_{0} \otimes id ) U_{(12)} ( e_{l \lambda i^{\prime}} \otimes 1 \otimes q_{l m \lambda i^{\prime}} ) ~~, ~~ e_{j \lambda i_{0}} \otimes 1 \ot 1 \right\rangle \eean

 \bean = \sum_{k, k^{\prime}}  \left\langle ( P^{\bot}_{0} \otimes id ) ( e_{k^{\prime} \lambda 0} \otimes q_{k^{\prime} k \lambda 0} \otimes q_{k  m \lambda 0} ) ~~, ~~ e_{j \lambda i_0} \otimes 1 \ot 1 \right\rangle \eean  
 
\bean  + \sum_{i^{\prime \prime}  \neq 0,~ k,~ k^{\prime}} \left\langle ( P^{\bot}_{0} \otimes id ) (   e_{k^{\prime \prime} k \lambda i^{\prime}} \otimes q_{k m \lambda 0} ) ~~ , ~~ e_{j \lambda i_0} \otimes 1 \ot 1 \right\rangle   \eean
 
 \bean + \sum_{i^{\prime} \neq 0,~ l,~l^{\prime}} \left\langle ( P^{\bot}_{0} \otimes id ) (    e_{l^{\prime} \lambda i^{\prime}} \otimes q_{l^{\prime} l \lambda i^{\prime}} \otimes q_{l m \lambda i^{\prime}} ) ~~ , ~~ e_{j \lambda i_0} \otimes 1 \ot 1 \right\rangle \eean
 
 \bean + \sum_{i^{\prime} ~ \neq 0,~ i^{\prime \prime} \neq i^{\prime},~ l,~ l^{\prime \prime}} \left\langle ( P^{\bot}_{0} \otimes id ) (   e_{l^{\prime \prime} \lambda i^{\prime \prime}} \otimes q_{l^{\prime \prime} l \lambda i^{\prime \prime}} \otimes q_{l m \lambda i^{\prime}} ) ~~ , ~~ e_{j \lambda i_0} \otimes 1 \ot 1 \right\rangle \eean
 
 \bean = \sum_{i^{\prime} \neq 0,~ k^{\prime},~ k^{\prime \prime}} \left\langle   e_{k^{\prime \prime} \lambda i^{\prime}} \otimes q_{k^{\prime \prime} k \lambda i^{\prime}} \otimes q_{k m \lambda 0} ~~,~~ e_{j \lambda i_0} \otimes 1  \ot 1 \right\rangle   \eean 
 
 \bean + \sum_{i^{\prime} \neq 0,~ l,~ l^{\prime}} \left\langle  e_{l^{\prime} \lambda i^{\prime}} \otimes q_{l^{\prime} l \lambda i^{\prime}} \otimes q_{l m \lambda i^{\prime}}  ~~,~~ e_{j \lambda i_0} \otimes 1  \ot 1 \right\rangle   \eean
 
 \bean + \sum_{i^{\prime} \neq 0,~ i^{\prime \prime} \neq i^{\prime},~ i^{\prime \prime} \neq 0, ~ l,~ l^{\prime \prime} }  \left\langle e_{l^{\prime \prime} \lambda i^{\prime \prime}} \otimes q_{l^{\prime \prime} l \lambda i^{\prime \prime}} \otimes q_{l m \lambda i^{\prime}} ~~,~~  e_{j \lambda i_0} \otimes 1  \ot 1 \right\rangle \eean
 
$ \in \cli \otimes QISO^{+}_{I}( d + d^* ) + QISO^{+}_{I}( d + d^* ) \otimes \cli $ ( as $ q_{k j \lambda i^{\prime}} \in \cli $ if $ i^{\prime} \neq 0  )$ \qed

%
%
%
%
%
%
%
%

\bthm

If $ \alpha_{U_{d + d^{*}}} $ is a $ C^* $ action on $ \cla ,$ then we have $ \clq^{\cll_{D}} \cong QISO^{+}_{I} ( d + d^* ) / \cli .$

\ethm

{\it Proof :}\\  By   Proposition \ref{Q_L_leq_QISO(d + d*)},
 we conclude that there exists a surjective CQG morphism $ \pi : QISO^{+}_{I} ( d + d^* ) \rightarrow \clq^{\cll_{D}}.$ 
By construction( Section 2.4 in \cite{goswami} ), the unitary representation  $ U_\cll $ of $ \clq^{\cll_{D}} $ preserves each of the $ \clh^{i}_{D},$ in particular
 $\clh^0_D$. It is then clear from the definition of $\cli$ that  
 $ \pi $ induces a surjective CQG morphism (in fact, a morphism in the category ${\bf Q}_I^\prime(d+d^*)$) $ \pi^{\prime} : QISO^{+}_{I} ( d + d^* ) / \cli \rightarrow \clq^{\cll_{D}}.$

Conversely,  if $ V=({\rm id} \ot \rho_\cli)\circ U_{d+d^*}$ is the representation of $QISO^{+}_{I} ( d + d^* )/\cli $ on $ \clh_{d + d^*} $ 
induced by $U_{d+d^*}$ (where  $\rho_\cli: QISO^+_I(d+d^*)\raro QISO^+_I(d+d^*)/\cli$ denotes the quotient map), then $V$ preserves $\clh^0_D$ (by definition of 
$\cli$), so commutes with $P_0$. Since $V$ also commutes with $(d+d^*)^2$, it follows that $V$ must commute with $(d+d^*)^2P_0=\cll$, i.e.
 $$ \widetilde{V} ( d^*d ~ P_0 \otimes 1 ) = ( d^*d P_0 \otimes 1 ) \widetilde{V}. $$  It 
 is easy to show from the above that $\alpha_V$ (which is a $C^*$ action 
   on $\cla$ since $\alpha_{U_{d+d^*}}$ is so by assumption) is a smooth isometric action of $QISO^+_I(d+d^*)/\cli$ in the sense of \cite{goswami}, w.r.t. 
 the Laplacian $\cll$. This implies that $QISO^+_I(d+d^*)/\cli$ is a sub-object of $\clq^\cll$ in the category ${\bf Q}^\cll$, and completes the proof.
\qed
\vspace{2mm}\\ 

%

Now we prove that under some further assumptions which are valid for classical manifolds as well as their Rieffel deformation,
  one even has the isomorphism $ \clq^{\cll_{D}} \cong QISO^{+}_{I}( d + d^* ).$

We assume the following:\\

({\bf A})  Both the spectral triples $ ( \cla^{\infty}, ~ \clh, ~ D ) $ and $ ( \cla^{\infty}, ~ \clh_{d + d^{*}}, ~ d + d^{*} ) $ satisfy the assumptions $({\bf 1})-({\bf 7})$, so in particular
 both  $\clq^{\cll_D}$ and $\clq^{\cll_{D^\prime}}$ exist (here $D^\prime=d+d^*$).\\
({\bf B}) $$ {\left\langle a,~ b \right\rangle}_{\clh^{0}_{ D }} = {\left\langle a,~ b \right\rangle}_{\clh^{0}_{ D^{\prime} }} $$
 and $$ {\left\langle d_{D}a,~ d_{D}b \right\rangle}_{\clh^{1}_{ D }} = {\left\langle d_{D^{\prime}}a,~ d_{D^{\prime}}b 
\right\rangle}_{\clh^{1}_{ D^{\prime} }} \forall a, ~ b \in \cla^{\infty}.$$

\brmrk

For classical compact spin  manifolds these assumptions can be verified by comparing the local expressions of $D^2$ and the `Hodge Laplacian' $(D^\prime)^2$ in suitable coordinate charts. In fact, in this case, both these operators turn out to be essentially same, upto a `first order term', which is relatively compact w.r.t. $D^2$  or $(D^\prime)^2$.


\ermrk

 From assumption $({\bf  B})$, we observe that the identity map on $ \cla^{\infty} $ extends to a unitary, say $ \Sigma $, from $ \clh^{0}_{D} $ to 
$ \clh^{0}_{D^{\prime}}.$ Moreover,  we have  $$ \cll_{D} = \Sigma^{*}_{0} \cll_{D^{\prime}} \Sigma,$$ from which the following follows immediately:

\bppsn

\label{Q_L_D_isom_Q_L_(d + d*)}

Under the above asumptions, $ \clq^{\cll_{D}} \cong \clq^{\cll_{D^{\prime}}} $

\eppsn


%
%
We conclude this section with the following result, which identifies the quantum isometry group $\clq^{\cll_D}$ of \cite{goswami} as the $QISO^+_I$ of a spectral triple, and thus, in some sense, acommodates the construction of \cite{goswami} in the framework of the present article. 
\bthm

\label{QISO( D ) vs  QISO( L )}

If in addition to the assumptions already made, $(\cla^\infty, \clh_{D^\prime}, D^\prime)$ also satisfies the conditions of  Theorem \ref{QISO_I_D_<_Q_L_D}, so that $QISO^+_I(D^\prime)$ has a $C^*$-action, the we have $$ \clq^{\cll_{D}} \cong QISO^{+}_{I}( D^{\prime} ) \cong \clq^{\cll_{D^{\prime}}}.$$

\ethm

{\it Proof :}\\ 
By  Proposition \ref{Q_L_leq_QISO(d + d*)}  that we have $ \clq^{\cll_{D}} $ is a sub-object of $ QISO^{+}_{I}( D^{\prime} ) $ 
in the category ${\bf Q}^\prime_{I}( D^{\prime} ) $. 
On the other hand, by Theorem \ref{QISO_I_D_<_Q_L_D} we have $ QISO^+_{I}( D^{\prime} ) $ is a sub object of $ \clq^{\cll_{D^{\prime}}} $ in the category $ {\bf Q}^{\cll_{D^\prime}}.$
Comining these facts with the colclusion of Proposition \ref{Q_L_D_isom_Q_L_(d + d*)}, we get the required isomorphism.\qed




 

\brmrk

The assumptions, and hence the conclusions, of this section are valid also for spectral triples obtained by Rieffel deformation of a classical spectral triple, to be discuused in details in section 5.

\ermrk

\section{Examples and computations}
\subsection{ Equivariant spectral triple on $SU_\mu(2)$}
Let $ \mu \in [ -1, 1 ].$
 Then $ SU_{\mu}( 2 ) $ is defined as the universal unital $ C^{*} $ algebra generated by $\alpha,~ \gamma $ such that : 
  \be \label{su2def1}  {\alpha}^{*} \alpha + {\gamma}^{*}\gamma = 1 \ee
  
  \be \label{su2def2}  \alpha {\alpha}^{*} + {\mu}^{2}{\gamma}{\gamma}^{*} = 1 \ee
  
  \be \label{su2def3} {\gamma}{\gamma}^{*} = {\gamma}^{*}\gamma \ee
  
  \be \label{su2def4} \mu \gamma \alpha = \alpha \gamma \ee
  
  \be \label{su2def5} \mu {\gamma}^{*} \alpha = \alpha {\gamma}^{*} \ee
   
   Let $ \clh = L^{2}( SU_{\mu}( 2 ) ) $ be the G.N.S space of $ SU_{\mu}( 2 ) $ with respect to the Haar state $ h .$
   
   For each $ n \in \{ 0, 1/2,1,..... \},$ there is a unique irreducible representation $ T^{n} $ of dimension $ 2n + 1.$ Denote by $ t^{n}_{ij} $ the $ ij $ th entry of $ T^{n}.$ They form an orthogonal basis of  $ \clh.$ Denote by $e^{n}_{ij} $ the normalized $ t^{n}_{ij} $ s so that $ \{ e^{n}_{ij} : n = 0, 1/2,1,......,i,j = -n, -n +1,.....n \}$ is an orthonormal basis.
 
  Consider the  spectral triple on $ SU_{\mu}( 2 ) $ constructed by Chakraborty and Pal ( \cite{partha} ) and also discussed thoroughly in \cite{con2} which is defined by $ ( \cla^{\infty}, \clh, D ) $ where $\cla^\infty$ is the linear span of $t^n_{ij}$ s,  and $ D $ is defined by : 
  \bean
    \lefteqn{ D ( e^{n}_{ij} ) }\\
     &=& ( 2n + 1 )e^{n}_{ij},~   n \neq i \\
     &=&  -( 2n + 1 ) e^{n}_{ij},~ n = i \eean
      
 Here, we have a cyclic separating vector $ 1_{SU_{\mu}( 2 )}$, and the corresponding faithful state is  the Haar state $h$.                 
   Therefore, an operator commuting with $ D $ ( equivalently with $ \widehat{D} $ ) must keep $ V^{l}_{i}: = {\rm Span} \{ t^{l}_{ij} : j = -l,........l \} $ invariant $ \forall $ fixed $ l $ and $ i $ where $ \widehat{D} $ is the operator as in Section 2.3.

 In the notation of Corollary \ref{unrestrictedcorollary}
, we have $ \cla^{\infty}_{0} = span \{ t^{l}_{i,j} : l = 0, 1/2,......... \}.=\cla^\infty$ in this case. 
All the conditions of Theorem \ref{unrestricted} and Corollary \ref{unrestrictedcorollary} are satisfied. Thus,  the universal object of the category $ \widehat{\bf C} $ exists ( notation as in Corollary \ref{unrestrictedcorollary} ) and we denote it by $ \widehat{\clq}.$

We recall from \cite{klimyk} that 

\be t^{1/2}_{-1/2, -1/2} = \alpha, t^{1/2}_{-1/2, 1/2} = - \mu {\gamma}^{*} , t^{1/2}_{1/2, -1/2} = \gamma, t^{1/2}_{1/2, 1/2} =  {\alpha}^{*} \ee

Moreover, if $ f_{n,i} = a_{n,i} \alpha^{n - i} \gamma^{* n + i} $ ( where $ a_{n,i} $ s are some constants as in \cite{klimyk} ) then  $ \{ f_{n,i} : n = 0, \frac{1}{2}, 1, \frac{3}{2},....., ~ -n \leq i \leq n \} $ is an orthonormal basis of $ SU_{\mu}( 2 ) .$

Now, $ f_{n + \frac{1}{2},i} = c ( n,i ) \alpha f_{n, i + \frac{1}{2}} $ for some constants $ c_{n,i} .$ Applying the coproduct on both sides and then comparing coefficients, we have the following recursive relations.

\bea 
\label{rec_1} \lefteqn{ t^{l + 1/2}_{i, l + 1/2} }{\nonumber}\\
 &=&  c_{11}(i,l) t^{l}_{i 1/2, l} {\gamma}^{*} + c_{12}(i,l) t^{l}_{i- 1/2, l}{\alpha}^{*} ,~~~~~~~~ -l + 1/2 \leq i \leq l - 1/2 {\nonumber}\\
&=&  c_{21}(i,l) {\gamma}^{*} t^{l}_{i+ 1/2, l}   ~~~~~~~~~~         i = -l - 1/2 {\nonumber}\\ 
&=&  c_{31}(i,l) {\alpha}^{*} t^{l}_{i - 1/2, l} ~~~~~~~~          i= l + 1/2  {\nonumber}\\
\eea

For $ j \le l,$

\bea \label{rec_2}  \lefteqn{ t^{l + 1/2}_{i, j}  }{\nonumber}\\
 & =&  c(l,i,j) \alpha t^{l}_{i+ 1/2, j + 1/2}  +c^\prime(l,i,j)  \gamma t^{l}_{i - 1/2, j + 1/2} ~~~~~~~    - l + 1/2 \leq i \leq l - 1/2 {\nonumber}\\
& = & d(l,j) \alpha  t^{l}_{- l, j + 1/2} + d^{\prime}(l,j) \gamma^* t^{l}_{- l, j - \frac{1}{2}} ~~~~~~~   i = -l - 1/2, ~ - l + \frac{1}{2} \leq j \leq l - \frac{1}{2} {\nonumber}\\
& = & d^{\prime \prime}(l,j) \alpha  t^{l}_{i+ 1/2, j + 1/2} ~~~~~~~   i = -l - 1/2, ~ j = - l - \frac{1}{2} {\nonumber}\\
& =&  e(l,j)  \gamma t^{l}_{i - 1/2, j+ 1/2} +  e^{\prime}( l,j ) \alpha^* t^{l}_{i - \frac{1}{2}, j - \frac{1}{2}} ~~~~~~ i= l + 1/2 {\nonumber};\\
\eea
where $C_{pq}(il), c(l,i,j), d(l,j), d^\prime_{l,j}, d^{\prime \prime}(l,j), e(l,j), e^{\prime}( l,j )  $ are all complex numbers.

\blmma

\label{equivalent}

Given a CQG $\clq$ with an action $\Phi$ on $\cla$, the following are equivalent :

1.$ ( \clq, \Phi ) \in Obj ( \widehat{{\bf C}} ) $

2. The action is linear, i.e, $ V^{1/2}_{-1/2} $ ( equivalently $,V^{1/2}_{1/2} $ ) is invariant under $ \Phi $ and the  representation obtained by restricting  $\Phi$  to $ V^{1/2}_{1/2} $ is a unitary representation.

3. $ \Phi $ is linear and Haar state preserving.

4. $ \Phi $ keeps $ V^{l}_{i} $ invariant for each fixed $ l $ and $ i .$

\elmma

{\it Proof :}\\  

$ 1. \Rightarrow 2. $ Since  $ \Phi $ commutes with $ \widehat{D} $, $ \Phi $ keeps each of the eigenspaces of $ \widehat{D} $ invariant and so in particular  preserves $ V^{1/2}_{-1/2} $ , i.e $ \Phi $ is linear. The condition that $ ( h \otimes id ) \Phi = h( \cdot).1 $ implies that the corresponding representation induces a unitary.

\vspace{4mm}

$ 2 \Rightarrow 3. $ By linearity, write  $ \Phi ( \alpha ) = \alpha \otimes X + {\gamma}^{*} \otimes Y $ and $ \Phi ( {\gamma}^{*} ) = \alpha \otimes Z + {\gamma}^{*} \otimes W.$

Firstly, $ \Phi $-invariance of ${\rm Span} \{t^{k}_{i,j}\} $ for $ k = 0 $ and $k= \frac{1}{2} $ follow from  linearity and the fact that $ \Phi( 1 ) = 1.$ 

Next, we show that $ \Phi $ keeps $ span \{ t^{1}_{ij} : i,j = -1,0,1 \} $ invariant.

We recall that $ t^{1}_{ij} $ is given by the matrix: 

$ \left ( \begin {array} {cccc}
   {{\alpha}^{*}}^{2} & - ( {\mu}^{2} + 1) \alpha^{*} \gamma & - \mu {\gamma}^{2}  \\ {\gamma}^{*} {\alpha}^{*} &  1 - ( {\mu}^{2} + 1) {\gamma}^{*} \gamma & \alpha \gamma \\ - \mu {{\gamma}^{*}}^{2} & - ( {\mu}^{2} + 1) {\gamma}^{*} \alpha & {\alpha}^{2}  \end {array} \right ) $ 

By inspection, we see that $ \Phi( V^{1}_{i} ) \subseteq V^{1}_{i} \otimes \clg $ for $ i = -1, 1.$

Hence, it is enough to check the $\Phi$-invariance   for $ \alpha \gamma $ and $ 1 - ( {\mu}^{2} + 1 ) {\gamma}^{*} \gamma .$

Comparing coefficients in $ \Phi ( \alpha \gamma ),$ we can see that it belongs to $V^1_0$  if and only if $ X{Z}^{*} + Y{W}^{*} = 0.$ Similarly, in the case of $ 1 - ( {\mu}^{2} + 1 ) {\gamma}^{*} \gamma,$ we have the condition $Z{Z}^{*} + W{W}^{*} = 1.$
But these conditions follow from the unitarity of the matrix $ \left(  \begin {array} {cccc}
     X^{*}   &  Z^{*}  \\ Y^{*} & W^{*} \end {array} \right ) ,$ 
        which is nothing but the matrix corresponding  to the restriction of $ \Phi $ to $V^{1/2}_{1/2} .$ 
 Thus, $ \Phi $ keeps $ {\rm Span} \{ t^{1}_{ij} : i,j = -1,0,1 \} $ invariant.

Moreover, by using the recursive relations ( \ref{rec_1} ), ( \ref{rec_2} ) and the following multiplication rule   ( see \cite{klimyk} ),
$$ t^{l}_{i,j} t^{1/2}_{{i}^{'},{j}^{'}} = \sum_{k = -\left| l - 1/2 \right|,..... l + 1/2 } c_{k}( l,i,j,{{i}^{'},{j}^{'}} ) t^{k}_{ i + {i}^{'},j + {j}^{'}} ,$$  ($c_k(l,i,j,i^\prime, j^\prime)$ are scalars) we can easily observe  that $ \forall l \geq 3/2, ~ \Phi ( V^{l + 1/2}_{i} )  \subseteq V^{l - 1/2}_{i} \oplus V^{l + 1/2}_{i}.$

Using thse observations, we conclude that for  $ t^{l}_{ij} $  that $ \Phi $ maps  ${\rm  Span} \{ t^{l}_{ij} : l \geq 1/2 \} $ into itself. 

So, in particular, $ {\rm Ker}( h) = span \{ V^{l}_{i} , i = -l,,,,l, l \geq 1/2 \} $ is invariant under $ \Phi $ which ( along with $ \Phi ( 1 ) = 1 $ ) implies that $ \Phi $ preserves $ h.$
 
\vspace{4mm}

$3. \Rightarrow 4.$  

 We proceed by induction.
 The induction hypothesis holds for $ l =  \frac{1}{2}$ since  linearity means that span $ \{ \alpha, {\gamma}^{*} \} $ is invariant under $ \Phi $ and hence span $ \{ {\alpha}^{*}, \gamma  \} $ is also invariant. The case for $ l = 1 $ can be checked by inspection as in the proof of $2 \Rightarrow 3$.
 Take the induction hypothesis that $\Phi$ keeps $V^k_i$ invariant  for all $k,i$ with $k \leq l$. From the proof of $ 2 \Rightarrow 3$ we also have $ \forall l \geq \frac{3}{2},$ $ \Phi (  V^{l + 1/2}_{i} )  \subseteq V^{l - 1/2}_{i} \oplus V^{l + 1/2}_{i},$ by using  linearity only. 
  Thus, $\widetilde{\Phi}$ leaves invariant the Hilbert $\clq$ module $(V^{l-\frac{1}{2}}_i \oplus V^{l+ \frac{1}{2}}_i) \ot \clq$, and is a unitary there since $\Phi$ is Haar-state preserving. Since $\widetilde{\Phi}$ leaves invariant $V^{l-\frac{1}{2}}_i \ot \clq$ by the induction hypothesis, it must keep its orthocomplement, i.e. $V^{l +\frac{1}{2}}_i$ invariant as well. 
  
  \vspace{4mm}

$ 4. \Rightarrow 3.$

 The fact that $ \Phi $ keeps $ V^{l}_{i} $ invariant for $ l = 1/2 $ will imply that $ \Phi $ is linear. The proof of haar state preservation is exactly the same as in $ 2 \Rightarrow 3.$
 
 \vspace{4mm}

$ 4 \Rightarrow 1.$ 

That  $ \Phi $ preserves the Haar state follows from arguments used before.  Since  $ \cla^{\infty}_{0} = 
 {\rm Span} \{ t^{l}_{ij} : l \geq 0, i,j = -l,.......l \},$ and  $ \Phi $ keeps $ V^{l}_{i} $ invariant it is obvious that $\Phi(\cla^\infty_0) \subseteq \cla^\infty_0 \ot_{\rm alg} \clq_0$ and   $ \Phi \widehat{D} = ( \widehat{D} \otimes id ) \Phi .$

 \qed

We now introduce the compact quantum group $U_{\mu}( 2 ).$ ( We refer to \cite{klimyk} for more details )

\bdfn

As a unital $ C^{*} $ algebra, $ U_{\mu}( 2 ) $ is generated by $4$ elements $ u_{11}, u_{12}, u_{21}, u_{22} $ such that:

$u_{11}u_{12} = \mu u_{12}u_{11}, u_{11}u_{21} = \mu u_{21}u_{11}, u_{12}u_{22} = \mu u_{22}u_{12}, u_{21}u_{22} = \mu u_{22}u_{21}, u_{12}u_{21} = u_{21}u_{12}, u_{11}u_{22} - u_{22}u_{11} = ( \mu - {\mu}^{-1} )u_{12}u_{21} $ and the condition that the matrix 

$ \left(  \begin {array} {cccc}
     u_{11}   & u_{12}  \\ u_{21} & u_{22} \end {array} \right ) $ is a unitary.
 
 The  CQG structure is given by $ \Delta ( u_{ij} ) = \sum_{k = 1,2} u_{ik} \otimes u_{kj} , \kappa ( u_{ij} ) = {u_{ji}}^{*}, \epsilon ( u_{ij} ) = \delta_{ij} .$   
 
 \edfn
 
 \brmrk

Let the quantum determinant $ D_{\mu} $ be defined by $ D_{\mu} = u_{11}u_{22} - \mu u_{12}u_{21} = u_{22}u_{11} - {\mu}^{-1} u_{12}u_{21} .$ Then, $ {D_{\mu}}^{*}D_{\mu} = D_{\mu}{D_{\mu}}^{*} = 1 .$ Moreover, $ D_{\mu} $ belongs to the centre of $ U_{\mu} ( 2 ) .$

\ermrk

 By  Lemma \ref{equivalent}, we have identified the category $ \widehat{\bf C} $ with the category  of CQG having actions on $ SU_{\mu} ( 2 ) $  satisfying 3. of Lemma \ref{equivalent}. Let the universal object of this category be denoted by $ ( \widehat{\clq}, \Gamma ).$

Then by linearity we can write:

$$ \Gamma ( \alpha ) = \alpha \otimes A + {\gamma}^{*} \otimes B $$

$$ \Gamma ( {\gamma}^{*} ) = \alpha \otimes C + {\gamma}^{*} \otimes D $$

Now we exploit the homomorphism condition of $ \Gamma $ to get relations between $ A, B, C, D .$

\blmma

\label{SU_mu_2_alpha* alpha + gamma* gamma = 1}

\be \label {1} {A}^{*}A + C{C}^{*} = 1 \ee
\be \label {2} {A}^{*}A + {\mu}^{2}C{C}^{*} = {B}^{*}B + D{D}^{*} \ee
\be \label {3} {A}^{*}B = - \mu D{C}^{*} \ee
\be \label {4} {B}^{*}A = - \mu C{D}^{*} \ee

\elmma

{\it Proof :}\\ The proof follows from the relation ( \ref{su2def1} ) by comparing coefficients of $ 1, {\gamma}^{*}\gamma ,{\alpha}^{*}{\gamma}^{*} $ and $ \alpha \gamma $ respectively. \qed

\blmma

\label{SU_mu_2_alpha alpha* + mu2 gamma gamma* = 1}

\be \label {5} A {A}^{*} + {\mu}^{2} C{C}^{*} = 1 \ee
\be \label {6} B {B}^{*} + {\mu}^{2}D{D}^{*} = {\mu}^{2}1 \ee
\be \label {7} B {A}^{*} = - {\mu}^{2} D{C}^{*} \ee

\elmma

{\it Proof :}\\ From the equation ( \ref{su2def2} ) by equating coefficients of $ 1 $ and $ {\alpha}^{*} {\gamma}^{*} ,$ we get respectively ( \ref{5} ) and ( \ref{7} ) whereas ( \ref{6} ) is obtained by equating coefficients of $ {\gamma}^{*} \gamma $ and using ( \ref{5} ).

  \qed

\blmma

\label{SU_mu_2_gamma* gamma = gamma gamma*}

\be \label {8}  {C}^{*} C = C {C}^{*} \ee
\be \label {9} ( 1 - {\mu}^{2} ) {C}^{*}C = {D}^{*}D - D{D}^{*} \ee
\be \label {10}  {C}^{*}D =  \mu D{C}^{*} \ee

\elmma

{\it Proof :}\\ The proof follows from the equation ( \ref{su2def3} ) from the coefficients of $ 1,{\gamma}^{*}\gamma, {\alpha}^{*} {\gamma}^{*}, $ respectively. \qed

\blmma

\label{SU_mu_2_alpha gamma = mu gamma alpha}

\be \label {11}  -{\mu}^{2} A{C}^{*} + B{D}^{*} - \mu D^{*}B + \mu C^{*} A = 0 \ee
\be \label {12}  A{C}^{*} = \mu {C}^{*}A  \ee
\be \label {13}  B{C}^{*} =  {C}^{*}B \ee 
\be \label {14}  A{D}^{*} = {D}^{*} A \ee

\elmma

{\it Proof :}\\ The proof follows from the equation ( \ref{su2def4} ) from the coefficients of $ {\gamma}^{*}\gamma, 1, {\alpha}^{*} {\gamma}^{*} $ and $ \alpha \gamma $ respectively. \qed

\blmma

\label{SU_mu_2_alpha gamma* = mu gamma* alpha}

\be \label {15}  AC = \mu CA  \ee
\be \label {16}  BD =  \mu DB \ee 
\be \label {17}  AD -\mu CB = DA -  {\mu}^{-1} BC = 0  \ee

\elmma

{\it Proof :}\\ The proof follows from ( \ref{su2def5}  ) from the coefficients of $ {\alpha}^{2}, {\gamma}^{*2}, {\gamma}^{*}\alpha $ respectively. \qed

\vspace{8mm}

Now we consider the antipode, say $ \kappa .$

From the condition $ ( h \otimes id )\Gamma( a ) = h( a ).1,$ we have that $ \Gamma $ gives a unitary representation of the compact quantum group via $ \widetilde{\Gamma} ( a \otimes q ) = \Gamma( a )( 1 \otimes q ) .$

Now, the restriction of this unitary representation to the orthonormal basis $ \{ \alpha, {\mu}^{-1}{\gamma}^{*} \} $ is given by the matrix : 
$ \left ( \begin {array} {cccc}
   A & \mu C  \\ {\mu}^{-1} B & D \end {array} \right ) .$ 
                        
 Similarly, with repect to the orthonormal set $ \{ {\alpha}^{*}, \gamma \}  $ ,this representation is given by the matrix: 
    $ \left ( \begin {array} {cccc}
   A^{*} & C^{*}  \\ B^{*} & D^{*} \end {array} \right ) .$

Thus, we have,
 
 $ \kappa( A ) = {A}^{*},~ \kappa( D ) = {D}^{*}, ~ \kappa( C ) = {\mu}^{-2} B^{*}, ~ \kappa( B ) = {\mu}^{2} C^{*}, ~ \kappa( A^{*} ) = A, ~ \kappa(C^{*} ) = B, ~ \kappa(B^{*}) = C, ~ \kappa( D^{*} ) = D .$
 
 Now, we apply $\kappa$ to the above equations to get the following additional relations.
 
 \blmma
 
 \be \label {27} AB = \mu BA \ee
 \be \label {28} CD = \mu DC \ee
\be \label {31}  B{C}^{*} = {C}^{*}B \ee

  \elmma

 {\it Proof :}\\ ( \ref{27} ), ( \ref{28} ), ( \ref{31} ) follow by applying $ \kappa $ to the equations ( \ref{15} ), ( \ref{16} ) and ( \ref{13} ) respectively.  
 
 
 \blmma
 
 \label{oneway}
 
 The map $ \phi : U_{\mu}( 2 ) \rightarrow \widehat{\clq} $ defined by $ \phi ( u_{11} ) = A, ~ \phi( u_{12} ) = \mu C, ~ \phi( u_{21} ) = {\mu}^{-1} B, ~ \phi( u_{22} ) = D $ is a $ \ast $ homomorphism.            

 \elmma
 
 {\it Proof :}\\ It is enough to check that the defining relations of  $ U_{\mu}( 2 ) $ are satisfied.
 
 1. $  \phi( u_{11} u_{12} ) = \phi( \mu u_{12} u_{11} ) \Leftrightarrow  \phi ( u_{11} ) \phi( u_{12} ) = \mu \phi( u_{12} ) \phi( u_{11} ) \Leftrightarrow  A ( \mu C ) = \mu ( \mu C )A \Leftrightarrow  AC = \mu CA $ which is the same as ( \ref{15} ).
 
 2. $  \phi( u_{11} u_{21} ) = \phi( \mu u_{21} u_{11} ) \Leftrightarrow A( {\mu}^{-1}B ) = \mu ( {\mu}^{-1} B ) A \Leftrightarrow AB = \mu BA $ which is the same as equation ( \ref{27} ).
 
 3. $  \phi ( u_{12} u_{22} ) = \phi( \mu u_{22} u_{12}  ) \Leftrightarrow \mu CD = \mu D ( \mu C ) \Leftrightarrow CD = \mu DC $ which is the same as equation ( \ref{28} ).
 
 4. $  \phi(u_{21} u_{22} ) = \phi( \mu u_{22} u_{21}  ) \Leftrightarrow {\mu}^{-1} BD = \mu D {\mu}^{-1} B \Leftrightarrow BD = \mu DB $ which is the same as equation ( \ref{16} ).

 5.$ \phi( u_{12}u_{21}  ) = \phi( u_{21}u_{12}    ) \Leftrightarrow \mu C {\mu}^{-1} B = {\mu}^{-1} B \mu C \Leftrightarrow CB = BC . $ 
  
  Now, $ B{C}^{*} = {C}^{*}B $ from equation ( \ref{31} ). But from ( \ref{8} ), $ C $ is normal. hence $ BC = CB .$

 6.  $  \phi( u_{11} u_{22} - u_{22} u_{11}) = ( \mu - {\mu}^{-1} )\phi( u_{12} u_{21}) \Leftrightarrow AD - DA = ( \mu - {\mu}^{-1} ) \mu C {\mu}^{-1} B.$
 
 From ( \ref{17} ), we have $ AD - DA = \mu CB - {\mu}^{-1} BC.$
 
 But $ B{C}^{*} = {C}^{*}B $ from ( \ref{31} ) and $ C $ is normal from ( \ref{8} ). Hence $ BC = CB.$ Hence $ AD - DA = ( \mu - {\mu}^{-1} ) CB $ holds.
 
 \qed

\blmma

 There is a $ C^{*} $ action $ \Psi $ of $ U_{\mu} ( 2 ) $ on $ SU_{\mu} ( 2 ) $ such that $ ( U_{\mu}( 2 ), \Psi ) \in Obj( \widehat{\bf C} ) .$

\elmma

{\it Proof :}\\ Define a $ \ast $ homomorphism $ \Psi $ on $ SU_{\mu}( 2 ) $ by 

    $$ \Psi( \alpha ) = \alpha \otimes u_{11} + {\gamma}^{*} \otimes \mu u_{21} $$
     
   $$  \Psi( {\gamma}^{*} ) = \alpha \otimes {\mu}^{-1}u_{12} + {\gamma}^{*} \otimes  u_{22} $$
   
   The homomorphism conditions are exactly the conditions ( \ref{1} ) - ( \ref{17} )  with $ A,B,C,D $ replaced by $ u_{11}, \mu u_{21}, {\mu}^{-1}u_{12} $ and $ u_{22} $ respectively.We check one of the relations and remark that the proof of the others are exactly similar.
    We prove ( \ref{1} )  i.e, $ {u_{11}}^{*}u_{11} + ( {\mu}^{-1}u_{12} ) { ({\mu}^{-1} u_{12} ) }^{*} = 1.$
    
  Using the fact that $ D_{\mu} $ is a central element of $ U_{\mu} ( 2 ) ,$ we have the left hand side $ =  u_{22} D^{-1}_{\mu} u_{11} + {\mu}^{-2} u_{12} ( - \mu u _{21} D^{-1}_{\mu} )  = ( u_{22} u_{11} - {\mu}^{-1} u_{12} u_{21} )D^{-1}_{\mu}  =  D_{\mu} D^{-1}_{\mu}  = 1 = $ right hand side.

 Clearly, $ \Psi $ keeps  span $ V^{1/2}_{-1/2} $ invariant and the corresponding representation is a unitary.
    
    Hence, by Lemma \ref{equivalent},  $ ( U_{\mu}( 2 ), \Psi ) \in Obj ( \widehat{\bf C} ) .$
    
    \qed
    
    \bcrlre
    
    There exists a surjective compact quantum  group morphism from $ \widehat{\clq} $ to $ U_{\mu}( 2 ) $ sending $ A, \mu C, {\mu}^{-1} B, $ and $ D $ to $ u_{11}, u_{12}, u_{21} $ and $ u_{22} $ respectively.  
    
    \ecrlre
    
    \bthm

  $  \widehat{Q} = U_{\mu}( 2 ) $ and hence $ \widetilde{{QISO}^{+}}( D ) = U_{\mu}( 2 ) \ast C( \IT ) $
    
    \ethm
   
   {\it Proof :}\\ Follows from  Lemma \ref{oneway} and the Corollary.
    The second part follows from Theorem \ref{unrestricted}.
    \qed

\subsection{The Podles Spheres}

Let $ SU_{\mu}( 2 ) $ be as in the previous section.

We recall the definition of the Podles sphere as in \cite{Dabrowski_et_al}.

We take $ \mu \in ( 0 , 1 )$ and $ t \in ( 0, 1 ].$
Let $[n] \equiv  [ n ]_{\mu} = \frac{\mu^n - \mu^{-n}}{\mu - \mu^{-1}} ), ~ n \in \IR. $

Let $ S^{2}_{\mu, c} $ be the universal $ C^{*} $ algebra generated by elements $ x_{-1}, x_{0}, x_{1} $ satisfying the relations:

\bean x_{-1} ( x_0 - t ) = \mu^2 ( x_0 - t ) x_{-1} \eean

\bean x_1 ( x_0 - t ) = \mu^{-2} ( x_0 - t ) x_1 \eean

\bean - [ 2 ] x_{-1} x_1 + ( \mu^2 x_0 + t ) ( x_0 - t ) = {[ 2 ]}^{2} ( 1 - t ) \eean

\bean - [ 2 ] x_1 x_{-1} + ( \mu^{-2} x_0 + t ) ( x_0 - t ) = {[ 2 ]}^{2} ( 1 - t ) \eean 

where $ c = t^{-1} - t, t > 0. $

The involution on $ S^{2}_{\mu, c} $ is given by 

\bean  x^{*}_{-1} = - {\mu}^{-1} x_1, ~~ x^{*}_{0} = x_0 \eean

We note that $ S^{2}_{\mu,c} $ as defined above is the same as $ \chi_{q, \alpha^{\prime}, \beta } $ of \cite{klimyk}( page 124 ) with $ q = \mu, ~ \alpha^{\prime} = t, ~ \beta = t^2 + \mu^{-2} {( \mu^2 + 1 )}^{2} ( 1 - t ) .$

Thus, from the expressions of $ x_{-1}, x_0, x_1 $ given in Page 125 of \cite{klimyk}, it follows that $ S^{2}_{\mu,c} $ can be realized as a $ \ast $ subalgebra of $ SU_{\mu}( 2 ) $ via:

\be \label{sphere_x_-1_in_termsof_su_mu2} x_{-1} = \frac{\mu {\alpha}^{2} + \rho ( 1 + \mu^2 ) \alpha \gamma - \mu^2 \gamma^2 }{ \mu ( 1 + \mu^2 )^{\frac{1}{2}} } \ee

\be \label{sphere_x_0_in_termsof_su_mu2} x_0 = - \mu \gamma^* \alpha + \rho ( 1 - ( 1 + \mu^2 ) \gamma^* \gamma ) - \gamma \alpha^* \ee

\be \label{sphere_x_1_in_termsof_su_mu2} x_1 = \frac{ \mu^2 {\gamma}^{*2} - \rho \mu ( 1 + \mu^2 ) \alpha^* \gamma^* - \mu \alpha^{*2} }{( 1 + \mu^2 )^{\frac{1}{2}} } \ee

where $ \rho^2 = \frac{\mu^2 t^2}{ ( \mu^2 + 1 )^2 ( 1 - t ) } .$

Setting

\bean A = \frac{ 1 - t^{-1} x_0}{ 1 + \mu^2}, ~~ B = \mu ( 1 + \mu^2 )^{- \frac{1}{2}} t^{-1} x_{-1} \eean

one obtains ( \cite{Dabrowski_et_al} ) that $ S^{2}_{\mu,c} $ is the same as the Podles' sphere as in \cite{podles}, i.e, the universal $ C^{*} $ algebra generated by elements $ A $ and $ B $ satisfying the relations:

$$ A^{*} = A,~  AB = \mu^{-2} BA,$$

$$ B^* B = A - A^2 + cI, ~ B B^* = \mu^2 A - \mu^4 A^2 + cI $$

We now introduce the spectral triple on $ S^{2}_{\mu c} $ as in \cite{Dabrowski_et_al}.

Let $ s = c^{\frac{1}{2}} \lambda^{- 1}_{+}, ~ \lambda_{+} = \frac{1}{2} + {( c + \frac{1}{4} )}^{\frac{1}{2}} $

$ \forall j \in \frac{1}{2} \IN,$

 $ u_j = ( {\alpha}^{*} - s \gamma^{*} ) ( \alpha^{*} - \mu^{-1} s \gamma^* )......( \alpha^* - \mu^{- 2j + 1} \gamma^* )$
 
 $ w_j = ( \alpha - \mu s \gamma ) ( \alpha - \mu^2 s \gamma )........( \alpha - \mu^{2j} s \gamma ) $
 
 $ u_{- j} = E^{2j} \triangleright w_j $
 
 $ u_0 = w_0 = 1 $
 
 $ y_1 = {( 1 + \mu^{- 2} )}^{\frac{1}{2}} ( c^{\frac{1}{2}} \mu^2 \gamma^{*2} - \mu \gamma^* \alpha^* - \mu c^{\frac{1}{2}} \alpha^{*2} ) $
 
 $ N^{l}_{kj} = {\left\| F^{l - k} \triangleright ( {y_1}^{l - \left| j \right|} u_j ) \right\|}^{- 1}  $
 
 Define $ v^{l}_{k,j} = N^{l}_{k,j} F^{l - k} \triangleright ( y^{ l - \left| j \right| }_{1} u_j ), ~ l \in \frac{1}{2} \IN_{0}, ~ j,k = - l,......l .$
 
 Let $ \clm_{N} $ be the Hilbert space with orthonormal basis $ \{ v^{l}_{m,N} : l = \left| N \right|, ~ \left| N \right| + 1,~ ........,~ m = - l,.......l \}.$
 
 Let
 
  $$ \clh = \clm_{- \frac{1}{2}} \oplus \clm_{\frac{1}{2}} .$$
 
 A representation $ \pi $ of $ S^{2}_{\mu,c} $ on $ \clh $ is defined by
 
$$ \pi ( x_i ) v^{l}_{m,N} = \alpha^{-}_{i} ( l, m; N ) v^{l - 1}_{m + i, N} + \alpha^{0}_{i} ( l, m; N ) v^{l}_{m + i, N} + \alpha^{+}_{i} ( l,m; N ) v^{l + 1}_{m + i,N} $$

where $ \alpha^{-}_{i},~ \alpha^{0}_{i}, ~ \alpha^{+}_{i} $ are as defined in \cite{Dabrowski_et_al}.

Finally as in Proposition 7.2 of \cite{Dabrowski_et_al}, the Dirac operator is defined by

$$ D ( v^{l}_{m, \pm \frac{1}{2}} ) = ( c_{1} l + c_2 ) v^{l}_{m, \mp \frac{1}{2}}  $$

where $ c_1, c_2 \in \IR, c_1 \neq 0 .$

Then $ ( S^{2}_{\mu,c}, \clh, D ) $ will be the spectral triple with which we are going to work.

In  \cite{sphere}, we have shown the following:

\bthm

Let the  positive, unbounded operator  $ R $ on $ \clh $ be defined by $ R ( v^{n}_{i,\pm \frac{1}{2}} ) = \mu^{- 2i} v^{n}_{i,\pm \frac{1}{2}} .$ Then,$ \tau_{R} $ equals $ h, $ i.e, the canonical Haar state on $ SU_{\mu} ( 2 ) $ and

 $$   {{QISO^{+}}_{R}} ( D ) = SO_{\mu} ( 3 ).  $$
 
 \ethm

 \brmrk

We have also worked with the spectral triple on the Podle's sphere $ S^{2}_{\mu 0} $ as in \cite{wagner}, and obtained the same result, i.e. identified $QISO^+_R$ of the spectral triple on $S^2_{\mu,0}$ with $SO_\mu(3)$.




   \ermrk

  We have also worked in \cite{sphere} with  another class of  spectral triple  introduced in \cite{chak_pal} for $c>0$.
  
Let $ \clh_{+} = \clh_{-} = l^2 ( \IN ),\clh = \clh_{+} \oplus \clh_{-} .$

Let $ e_n $ be an orthonormal basis of $ \clh_{+} = \clh_{-} $ and $ N $ be the operator defined on it by $ N ( e_n ) = n e_n .$

We recall the irreducible  representations $ \pi_{+} $ and $ \pi_{-}  : \clh_{\pm} \rightarrow \clh_{\pm} $   as in \cite{chak_pal}.

$$ \pi_{\pm} ( A ) = \lambda_{\pm} {\mu}^2n e_n $$

$$ \pi_{\pm} ( B ) = {c_{\pm}( n )}^{\frac{1}{2}} e_{n - 1}  $$

where $ \lambda_{\pm} =  \frac{1}{2} \pm {( c + \frac{1}{4} )}^{\frac{1}{2}}, c_{\pm} ( n ) = \lambda_{\pm} \mu^{2n} - {( \lambda_{\pm} \mu^{2n} )}^{2} + c .$

Let $ \pi = \pi_{+} \oplus \pi_{-} $ and $ D = \left ( \begin {array} {cccc}
   0 & N  \\ N & 0 \end {array} \right ) .$
   
Then $ ( S^{2}_{\mu,c}, \pi, \clh ) $ is a spectral triple.

We showed in \cite{sphere} that 

\bthm

 $ QISO^{+} $ of this spectral triple exists and  is isomorphic with $ C ( \IZ_{2} ) \ast C ( \IT )^{* \infty},$ where $C(\IT)^{* \infty}$ denotes the free product of countably infinitely many copies of $C(\IT)$.

 \ethm
\brmrk
The above example shows that unlike the classical case, where isometry groups are Lie groups and hence have faithful imbedding into a matrix group, $QISO^+$ in general may fail to be a compact matrix quantum group. In fact, it will be quite interesting to find conditions under which $QISO^+$ will be so.
\ermrk

\subsection{A commutative example : spectral triple on $\IT^2$}

We consider the spectral triple $ ( \cla^{\infty}, \clh, D ) $ on $ \IT^{2} $ given by $ \cla^{\infty} = C^{\infty} ( \IT^{2} ), \clh = L^{2} ( \IT^{2} ) \oplus L^{2} ( \IT^{2} )$ and $ D =  
  \left ( \begin {array} {cccc}
    0 &  d_{1} + i d_{2}   \\  d_{1} - i d_{2} &  0  \end {array} \right ) ,$ 
    
    where we view $ C( \IT^{2} ) $ as the universal $ C^{*} $ algebra generated by two commuting unitaries $ U  $ and $ V, d_{1} $ and $ d_{2} $ are derivations on $ \cla^{\infty} $ defined by :
    
     \be d_{1} ( U ) = U, d_{1} ( V ) = 0, d_{2} ( U ) = 0, d_{2} ( V ) = V .\ee
    
 $ e_{1} = ( 1, ~ 0 ) $ and $ e_{2} = ( 0, ~ 1 ) $ form an orthonormal basis of the eigenspace corresponding to the eigenvalue zero.
  
  The Laplacian in the sense of \cite{goswami} is given by $ \cll ( U^{m}V^{n} ) = - ( m^{2} + n^{2} ) U^{m}V^{n}  .$  We recall that we denote the quantum isometry group from the Laplacian $ \cll $ in the sense of \cite{goswami} by $ Q^{\cll_{D}} .$
  
  \blmma
  
  Let $ W $ be a unitary  representation of a CQG $ \widetilde{\clq} $ which commutes with $ D.$ Then the induced action on $ C^{\infty}( \IT^{2} ) ,$ say  $ \alpha ,$ satisfies :
  
  \be \label{torus1} \alpha ( U ) = U \otimes z_{1} \ee
  
  \be  \label{torus2} \alpha ( V ) = V \otimes z_{2} \ee
   
   where $ z_{1}, z_{2} $ are two commuting unitaries.
  
  \elmma
 
 {\it Proof:}\\ We denote the the maximal Woronowicz $ C^{*} $ subalgebra of $ \widetilde{\clq} $ which acts on $ C ( \IT^{2} ) $ faithfully by $ \clq .$
 
 We observe that $ D^{2} \left (  \begin {array} {cccc}  a e_{1} \\ a e_{2} \end {array} \right ) =  \left ( \begin {array} {cccc}
    \cll ( a ) &  0   \\  0 &  \cll ( a )   \end {array} \right ) \left (  \begin {array} {cccc}   e_{1} \\  e_{2} \end {array} \right )  .$ 
    
 The fact that  $ U $ commutes with $ D $ implies that $ U $ commutes with $ {D}^{2} $ as well, and hence $ \alpha $ commutes with the Laplacian $ \cll.$
    Therefore, $ \clq $ is a quantum subgroup of $ \clq^{\cll_{D}}.$ From \cite{jyotish}, we conclude that $ \clq^{\cll_{D}} = C( \IT^{2}>\!\!\! \lhd {Z_{2}}^{3} ) $ Thus $ \clq $ must be of the form $C(G)$ for  a classical subgroup of the orientation preserving isometry group of  $ \IT^{2} $, which is easily seen to be $ \IT^{2} $ itself  and whose (co )action is given by $ U \mapsto U \otimes {z_{1}}^{\prime} $ and  $ V \mapsto V \otimes {z_{2}}^{\prime} $ where $ {z_{1}}^{\prime}, {z_{2}}^{\prime} $ are two unitaries  generating $ C( \IT^{2} ) .$ Hence the lemma follows.\qed
    
    \bthm
    
  $  {\widetilde{QISO}}^+( C^\infty( \IT^{2} ), \clh, D ) $ exists and is isomorphic with $ C( \IT^{2} ) \ast C ( \IT ) \cong C^*({\IZ}^2 \ast \IZ)$ (as a CQG). Moreover $QISO^+$ of this spectral triple is $C(\IT^2)$.
    
    \ethm
    
 {\it Proof:}\\  Let $ W $ be as in the previous Lemma. Then we have,
  
    \be \label{torus3} W ( e_{1} ) = e_{1} \otimes q_{11} + e_{2} \otimes q_{12} \ee
    
    \be \label{torus4} W ( e_{2} ) = e_{1} \otimes q_{21} + e_{2} \otimes q_{22} \ee
    
    By comparing coefficients of $ U e_{1} $ and $ U e_{2},$ in the both sides of the equality $ ( D \otimes id ) W ( U e_{1} ) =  W D U e_{1} ,$ we have,
    
   \be z_{1} q_{12} = z_{1} q_{21} \ee
   
   \be z_{1} q_{11} = z_{1} q_{22} \ee
   
  $ z_{1} $ is a unitary implies that $ q_{11} = q_{22} $ and $ q_{12} = q_{21}.$
   
   Similarly, from the relation $ ( D \otimes id ) W ( V e_{1} ) = W D V e_{1},$ we have $ q_{12} = - q_{21}, q_{22} = q_{11}.$
   
   By the above two sets of relations, we have :
   
  $ q_{12} = q_{21} = 0, q_{11} = q_{22} = q $ ( say )
  
  But the matrix $ \left ( \begin {array} {cccc}
    q_{11} &  q_{12}   \\  q_{21} &  q_{22}  \end {array} \right ) $
    
    is a unitary in $  M_{2} ( \widetilde{\clq} ) ,$ so  $ q  $ is a unitary.
    
    Moreover, we note that $ W ( a e_{i} ) = \alpha( a ) W ( e_{i} ) ~ \forall a \in C^{\infty}( \IT^{2} ) .$
    Using the previous Lemma and the above observations, we deduce that any CQG which has a unitary representation commuting with the Dirac operator is a quantum subgroup of $ C( \IT^{2} ) \ast C( \IT ).$
    
    Moreover, $ C( \IT^{2} ) \ast C( \IT ) $ has a unitary representation commuting with $ D ,$ given by the formulae ( \ref{torus1} ) - ( \ref{torus4} ) taking $ q_{12} = q_{21} = 0, q_{11} = q_{22} = q^{'} $ where $ q^{\prime} $ is the generator of $ C( \IT ) $ and  $ z_{1}, z_{2} $ to be the generator of $ C( \IT^{2} ).$ 
  This completes the proof. \qed
  
  \brmrk 
  
  The canonical grading on $ C ( \IT^2 ) $ is given by the operator $ ( id \otimes \gamma ) $ on $ L^{2} ( \IT^2 \otimes \IC^2 ) $ where $ \gamma $ is the matrix  $ \left ( \begin {array} {cccc}
    0 &  1  \\  - 1 &  0  \end {array} \right ) .$
    
    The representation of $ C( \IT^{2} ) \ast C( \IT ) $ clearly commutes with the grading operator and hence is isomorphic with  $ \widetilde{QISO} ( C ( \IT^2 ),  L^{2} ( \IT^2 \otimes \IC^2 ), D, \gamma ) .$
    
    \ermrk

 \brmrk
 
 This example shows that  the conditions of Theorem \ref{unrestricted} are not necessary for the existence of $  \widetilde{{QISO}^{+}}.$
 
 \ermrk

\subsection{Another class of commutative examples : the spheres}
We consider the usual Dirac operator on the classical $n$-sphere $S^n$. In fact, we shall first consider a slightly more general set-up as in the section 3.5, p.82-89 of \cite{friedrich}, which we very briefly recall here. Let $G$ be a compact Lie group, $K$ a closed subgroup, and let $M$ be the  homogeneous space $G/K$ with a  $G$-invariant metric. The algebra  $C^\infty(M)$ is identified with the algebra (say $\cla^\infty$) of $K$-invariant functions in $C^\infty(G)$,  i.e. functions $f$ satisfying $f(gk)=f(g)$ for all $g \in G, k \in K$.  The Lie algebra $\mathsf{g}$ of $G$ splits as a vector space direct sum $\mathsf{ g}=\mathsf{k} +\mathsf{m}$ where $\mathsf{k}$ is the Lie algebra of $K$ and $\mathsf{m}$ is suitable ${\rm Ad}(K)$-invariant subspace of $\mathsf{g}$ (see \cite{friedrich} for more details). Thus we have  the representation  of $K$ given by ${\rm Ad}: K \raro {\rm SO}(\mathsf{m})$, and the corresponding lift $\widetilde{{\rm Ad}}: K \raro {\rm Spin}(\mathsf{m})$.  The space of smooth spinors can then  be identified with the space of smooth  functions $\psi : G \raro \Delta $ satisfying $\psi(gk)= \kappa \widetilde{{\rm Ad}}(k^{-1})\psi(g)$ for all $g \in G,~ k \in K$, where $\kappa : {\rm Spin}(\mathsf{ m}) \raro {\rm GL}(\Delta)$ denotes the spin representation. The action of $C^\infty(M)$, identified with the $K$-invariant smooth functions on $G$,  is given by multiplication, and the  Dirac operator $D$ is given by $$ D \psi=\sum_{i=1}^m X_i \cdot  X_i(\psi),$$ where $m={\rm dim}(\mathsf{ m})$ and $\{ X_1,..., X_m \}$ is an orthonormal basis of $\mathsf{ m}$ with respect to the suitable invariant  inner product described in \cite{friedrich} and $\cdot$ denotes the Clifford multiplication. From this expression of $D$, it follows by using the fact that $X_i$ 's are acting as derivations on the algebra of smooth functions, that $[D,f]\psi=\omega_f \cdot \psi,$ where $\omega_f=\sum_i (X_i f) X_i$. In fact, the space $\Omega^1_D$, which is isomorphic with the (complexified) space of smooth $1$-forms on $M$,  can now be  identified with the space of (smooth) ${\rm Ad}_K$ invariant functions from $G$ to $\mathsf{m} \cong \IC^m$ (which is also isomorphic with $C^\infty(M) \ot \IC^m$), and the map $d_D$ is, w.r.t. this identification, is nothing but the map which sends $f \in \cla^\infty \cong C^\infty(M)$ to $\sum_i X_i(f) \otimes  X_i \in \cla^\infty \otimes \IC^m$. The Hilbert space of $1$-forms is isomorphic with $L^2(M) \otimes \IC^m$, where $L^2(M)$ is the Hilbert space  completion w.r.t. the $G$-invariant volume measure, and from the $G$-invariance of the volume measure it is clear that the adjoint $X_i^*$ of the (left invariant ) vector field $X_i$ (viewed as a closable unbounded map on $L^2(M)$) is $-X_i$. It follows that the Laplacian is given by,  $d_D^*d_D=-\sum_{i=1}^m X^*_i X_i$ on $\cla^\infty \cong C^\infty(M)$. It is in fact nothing but the Casimir $\Omega_G$ in the notation of \cite{friedrich}, since $Yf=0$ for any $Y \in \mathsf{k}$.

Now, we want to apply the above observations to the special case of $n$-spheres. The Laplacian on such spheres considered in \cite{Varilly_book} ( Page 17 ) is indeed the Casimir operator and so by Theorem 2.2 and Remark 3.3 of \cite{jyotish}  the corresponding quantum isometry group $QISO^\cll$ is commutative as a $C^*$ algebra. However, by Corollary \ref{QISO_I_D_<_Q_L} of the present paper, any object of the category ${\bf Q}^\prime_I$ must be a quantum subgroup of $ Q^{\cll_{D}} $, so is in particular commutative as a $C^*$ algebra, so must be of the form $C(G)$ for a subgroup $G$ of the universal group of orientation preserving (classical) Riemannian isometries of $S^n$, i.e. $SO(n+1)$. 
 To summarise, we have the following:
\bthm
The quantum isometry group ${QISO^+}_I(S^n)$ is isomorphic with $C(SO(n+1))$, i.e. coincides with the corresponding classical group.
\ethm
 \section{$QISO^+$ of deformed spectral triples}
    In this section, we give a general scheme for computing orientation-preserving  quantum isometry groups by proving that $\widetilde{{QISO}^{+}_{R}}$ of a deformed noncommutative manifold coincides with (under reasonable assumptions) a similar deformation of the $\widetilde{{QISO}^{+}_{R}}$ of the original manifold. The technique is very similar to the analogous result for the quantum isometry groups in terms of Laplacian discussed in \cite{jyotish}, so we often merely sketch the arguements and refer to a similar theorem or lemma in \cite{jyotish}. 
    
 We recall the generalities on compact quantum groups from Section 2.2. In particular, given a compact quantum group $(\cls,\Delta)$, the dense unital $\ast$-subalgebra $\cls_{0}$ of $\cls$ generated by the matrix coefficients of the irreducible unitary representations has a canonical Hopf $\ast$-algebra structure.  Moreover,  given an action $\gamma : \clb \raro \clb \ot \cls$ of the compact quantum group $(\cls, \Delta)$  on a unital $C^*$-algebra $\clb$, it is known that one can find a dense, unital  $\ast$-subalgebra $\clb_0$ of $\clb$ on which the action becomes an action by the Hopf $\ast$-algebra $\cls_{0}$. We shall use the Sweedler convention of abbreviating $\gamma(b) \in \clb_0 \ot_{\rm alg} \cls_{0}$ by $b_{(1)} \ot b_{(2)}$, for $b \in \clb_0$. This applies in particular to the canonical action of the quantum group $\cls$ on itself, by taking $\gamma=\Delta$. Moreover, for a linear functional $f$ on $\cls$ and an element $c \in \cls$ we shall define the `convolution' maps $f \diamond c :=(f \otimes {\rm id} ) \Delta ( c )$ and $c \diamond f := ({\rm id} \otimes f) \Delta ( c)$. We also define convolution of two functionals $f$ and $g$ by $(f \diamond g)(c)=(f \ot g)(\Delta(c))$.   
  We also need the following:    
       \bdfn Let $ (\cls, \Delta_\cls) $ be a compact quantum group.  
   A vector space $ M $  is said to be an algebraic  $ \cls $ co-module (or $\cls$ co-module) if  there exists a linear map $ \widetilde{\alpha} : M \rightarrow M \otimes_{\rm alg} \cls_0  $ such that\\
 1. $ ( \widetilde{\alpha} \otimes {\rm id} ) \widetilde{\alpha} = ({\rm  id}  \otimes \Delta_{\cls} ) \widetilde{\alpha} $\\
  2. $ ( {\rm id} \otimes \epsilon ) \widetilde{\alpha}(m) = \epsilon(m ) 1_{\cls} $ $\forall m \in M$.
  
  \edfn

\vspace{4mm}

  Let $ ( \cla, \IT^{n}, \beta ) $ be a $ C^{*} $ dyanamical system and $ \pi_{0} : \cla \rightarrow \clb( \clh ) $ be a faithful representation, where $ \clh $ is a separable Hilbert space.
  
   Let $ \cla^{\infty}  $ be the smooth algebra corresponding to the $ \IT^{n} $ action $ \beta.$ Then
 for each skew-symmetric $n \times n$ real matrix $J$, we refer to \cite{rieffel} for the construction of 
  the `deformed' $C^*$-algebra $\cla_J$ and their properties.

  Assume now that we are given  a spectral triple $ ( \cla^{\infty}, \pi_{0}, \clh , D ) $ of compact type. Suppose that $ D $ has eigenvalues $ \{ \lambda_{0}, \lambda_{1},......... \} $ and $ V_{i} $  denotes the (finite dimensional) eigenspace of $ \lambda_{i} $ and let $ \cls_{00} $ denote the linear span of $ \{ V_{i}: ~ i = 0, 1,2,.. \}$
  
 Suppose, furthermore,  that there exists a compact abelian lie group $ \widetilde{\IT^{n}} $, with a  covering map $ \gamma :  \widetilde{\IT^{n}} \rightarrow \IT^{n}  .$ The Lie algebra of both $ \IT^{n} $ and  $ \widetilde{\IT^{n}} $ are isomorphic with $ \IR^{n} $ and we denote by $ e $ and $ \widetilde{e} $ respectively the corresponding exponential maps, so that $ e( u ) = e( 2 \pi iu ) , u \in \IR^{n} $ and $ \gamma ( \widetilde{e}( u ) ) = e ( u ).$ By a slight abuse of notation we shall denote the $\IR^n$-action $\beta_{e(u)}$ by $\beta_u$. 
  
  We also make the following assumption:\\
   There exists a strongly continuous unitary representation $   V_{\tilde{g}}$ , $\tilde{g} \in \widetilde{\IT^{n}} $ of $  \widetilde{\IT^{n}} $ on $\clh$ such that \\
     (a) $  V_{\tilde{g}} D = D V_{\tilde{g}}$ $\forall \tilde{g}$,\\
   (b) $ V_{\tilde{g}} \pi_{0} ( a ) {V_{\tilde{g}}}^{-1} = \pi_{0} ( \beta_{g} ( a ) )$, where $a \in \cla, \tilde{g} \in   \widetilde{\IT^{n}},$ and $ g = \gamma ( \tilde{g} )  $.\\
   
 We shall now show that we can `deform' the given spectral triple along the lines of \cite{connes_etal}.   For each $J$, the map $ \pi_{J}: \cla^\infty \raro {\rm Lin}(\clh^\infty)$ (where $ \clh^{\infty} $ is the smooth subspace corresponding to the representation $ V $ and ${\rm Lin}(\clv)$ denotes the space of linear maps on a vector space $\clv$) defined by $$ \pi_{J} ( a ) s \equiv  a \times_{J} s :=\int \int \beta_{Ju} ( a ) \widetilde{\beta}_{v} ( s ) e( u.v ) du dv$$  extends to a faithful $\ast$-representation of the $C^*$-algebra $ \cla^{\infty} $ in $\clb( \clh )$  where   $ \widetilde{\beta_{v}} = V_{\widetilde{e}( v )} $ ( which clearly maps $ \clh^{\infty}$ into $\clh^{\infty} ).$  
    
    We can extend the action of $ \IT^n $ on the $ C^* $ subalgebra $ \cla_1 $ of $ \clb ( \clh ) $  generated by $ \pi_{0} ( \cla ), ~ \{ e^{itD} : t \in \IR \} $ and elements of the form  $ \{ [ D, a ]: a \in \cla^{\infty} \} $  by $  \beta_g ( X ) = V_{\widetilde{g}} X {V_{\widetilde{g}}}^{- 1} ~ \forall X \in ~ \cla_1 $ where by an abuse of notation, we denote the action by the same symbol $ \beta .$ Let $ {\cla_1}^{\infty} $  denote the smooth vectors of $ \cla $ with respect to this action.We note that $ \forall ~ a \in \cla^{\infty}_{1}, [ D, a ] \in {\cla}^{\infty}_{1} .$
    
    \blmma
    
    \label{qorient_deformation_extension}
    
    $ \beta $ is a strongly continuous action (in the $C^*$-sense)  of $ \IT^n $ on $ \cla_1 $ and hence   $ \forall ~ X \in { \cla_1 }^{\infty}, ~ \pi_{J} ( X ) $ defined by $ \pi_{J} ( X ) s = \int \int \beta_{Ju} ( X ) \widetilde{\beta}_{v} ( s ) e( u.v ) du dv $ is a bounded operator.
    
    \elmma
    
  {\it Proof:}\\ We note that $ \beta $ is already strongly continuous on the $ C^* $ algebra generated by $ \pi_{0} ( \cla ), ~ \{ e^{itD} : t \in \IR \} .$ Thus it suffices to check the statement for elements of the form  $ [ D, a ] $ where $ a \in \cla^{\infty} .$
  
  To this end, fix any one parameter subgroup $g_t$ of $G$ such that $g_t$ goes to the identity of $G$ as $t \raro 0$. 
Let $ T^{\prime}_{t}, ~ \widetilde{T_t} $ denote the group of normal $\ast$-automorphisms  on $ \clb(\clh) $ defined by $ T^{\prime}_{t} ( X ) = V_{g_{t}} X V_{g^{- 1}_{t}} $ and $ \widetilde{T_t} ( X ) = e^{i t D} X e^{- i t D} .$  As  $ V_{g_{t}} $ and $ D $ commute, so do 
their generators. In particular, each of these generators leave the domain of the other invariant. Note also that $\cla^\infty$ is in the domain of the both the generators, and the generator  of $\widetilde{T_t}$ is given by $[D, \cdot]$ there. Thus, for 
 $a \in \cla^\infty$, we have $a, [D, a] \in {\rm Dom}(\Xi)$ (where $\Xi$ is the generator of $T^\prime_t$),  and $\Xi([D,a])=[D,\Xi(a)] \in \clb(\clh)$.

Using this, we obtain 
  $ \left\| T^{\prime}_{t} ([ D, a ])  - [ D, a ] \right\|  = \int^{t}_{0}  T^{\prime}_{s} (\Xi( [ D, a ])) ds  \leq ~ t  \left\| \Xi([ D, a ]) \right\|  .$ 
 The required strong continuity follows from this. Then applying Theorem 4.6 of \cite{rieffel} to the $ C^{*} $ algebra $ \cla_1 $ and the action $\beta, $ we deduce that $ \pi_{J} ( X ) $ is a bounded operator. \qed
   
  \blmma
   
  For each $J$, $ ( \cla^{\infty}_{J}, \pi_{J}, \clh, D ) $ is a spectral triple, i.e, $ [ D, \pi_{J}( a ) ] \in \clb ( \clh ) $ $\forall ~ a  \in \cla^\infty_J $.
   
  \elmma 
  
  {\it Proof:}\\ $ [ D, \pi_{J} ( a ) ] ( s )
   = D \int \int \beta_{Ju} ( a ) \widetilde{\beta_{v}} ( s) e ( u. v ) du dv - \int \int \beta_{Ju} ( a ) \widetilde{\beta_{v}} ( D s ) e ( u. v ) du dv .$
   
   Using the expression  $ \int \int f ( u ) g ( v ) e ( u. v ) = lim \sum_{L} ( f ~ \psi_{m} )( p ) \widehat{g} ( p ) $ ( where notations are as in \cite{rieffel}, Page 20 ) and closability of $ D ,$ we have
   
    $ D \int \int \beta_{Ju} ( a ) \widetilde{\beta_{v}} ( s) e ( u. v ) du dv = \int \int D ( \beta_{Ju} ( a ) \widetilde{\beta_{v}} ( s) ) e ( u. v ) du dv .$
   
 Thus, the above expression equals $ \int \int D ( \beta_{Ju} ( a ) \widetilde{\beta_{v}} ( s) ) e ( u. v ) du dv $
 
 $ - \int \int \beta_{Ju} ( a ) D  \widetilde{\beta_{v}} ( s ) e ( u. v ) du dv $ 
  as $ D $ commutes with $ V. $ 
  

   So we have 
   
  $[D, \pi_J(a)](s) = \int \int [ D, \beta_{Ju} ( a) ] \widetilde{\beta_{v}} ( s) e ( u. v ) du dv
   = \int \int V_{\widetilde{Ju}} [ D, a ] {V_{\widetilde{Ju}}}^{- 1} \widetilde{\beta_{v}} ( s) e ( u. v ) du dv
   = \pi_{J} ( [ D, a ] ) $
   
   which is a bounded operator by Lemma \ref{qorient_deformation_extension}. \qed

   \blmma
   \label{1111}
   Suppose that $ ( \widetilde{\clq}, U ) ~ \in ~ {\rm Obj} ( {\bf Q}( \cla, \clh, D ) ) , $ and  
       there exists a unital $ \ast $-subalgebra $ \cla_{0} \subseteq \cla $ which is norm dense in every $ \cla_{J}$ such that \\
      $ \alpha_{U} ( \pi_{0} ( \cla_{0} )  ) \subseteq \pi_{0} ( \cla_{0} ) \otimes_{alg} \clq_{0}, $ where $\clq \subseteq \widetilde{\clq} $ is the smallest Woronowicz $ C^{*} $ subalgebra such that $ \alpha( \cla_{0} ) \subseteq \pi_{0} ( \cla_{0} ) \otimes \clq,$ and $ \clq_{0} $ is the Hopf $\ast$-algebra obtained by matrix coefficients of irreducible unitary (co)-representations of $\clq.$ Also, let $ S_0 = {\rm span} \{ a s : a \in \cla_0, s \in S_{00} \} ,$
Then we have the following:\\
(a) $U(\cls_0) \subseteq \cls_0 \ot_{\rm alg} \tilde{\clq}_0$.\\
(b) $\tilde{\alpha}:=U|_{\cls_0}: \cls_0 \raro \cls_0 \ot_{\rm alg} \tilde{\clq}_0$  makes $\cls_0$ an algebraic $\tilde{\clq}_0$ co-module, satisfying $$ \tilde{\alpha}(\pi_0(a) s)=\alpha_{U}(a) \tilde{\alpha}(s)~~\forall a \in \cla_0, ~s \in \cls_0.$$ 

Moreover, if  $ ~ C ( \widetilde{\IT^n} ) $ is a sub object of $ \widetilde{\clq} $ in $ {\bf Q}( \cla, \clh, D ) $, then $ C ( \IT^n ) $ is a quantum subgroup of $ \clq .$ 
\elmma
 {\it Proof:}\\  $ U $ commutes with $ D $ and hence preserves the eigenspaces of $ D $ which shows that $ U $ preserves $ S_{00}.$ Then, $ U ( a s ) = \alpha ( a ) U ( s ) \subseteq ( \cla_0 \otimes \clq_0 ) ( S_{00} \otimes \clq_0 ) \subseteq S_0 \otimes \clq_0.$ Thus, the first assertion follows.
 
  The second assertion follows from the definition of $ \widetilde{\alpha} $ and $ \alpha_{u}.$
 
We now prove the third assertion. Let us denote by $\gamma^*$ the dual map of $\gamma$, so that $\gamma^* : C(\IT^n) \raro C(\tilde{\IT}^n)$ is an injective $C^*$-homomorphism. It is quite clear that $({\rm id} \ot \pi_{\tilde{\clq}})  \circ \alpha (\cla_0) \subseteq {\rm Im}({\rm id} \ot \gamma^*)$, hence it follows that $\pi_{\tilde{\clq}}(\clq_0) \subseteq {\rm Im}(\gamma^*)$. Thus, $\pi_\clq:=(\gamma^*)^{-1} \circ \pi_{\tilde{\clq}} $ is a surjective CQG morphism from $\clq$ to $C(\IT^n)$, which identifies $C(\IT^n)$ as a quantum subgroup of $\clq$. 

\qed

\brmrk

 From the definitions of $ \cla_0 $ and $ S_0 ,$ it follows that\\
   (i)  $ \pi_{0} ( \cla_{0} ) \cls_{0} \subseteq \cls_{0} $,\\
   (ii) $ \beta_{g}( \cla_{0} ) \subseteq \cla_{0} ~ \forall g.$\\
   
 \ermrk  

\vspace{4mm}

Let us now fix the object $(\tilde{\clq}, U)$ as in the statement of Lemma \ref{1111}. 
 From now on, we will identify $ \cla^{\infty}_{J} $ with $ \pi_{J}( \cla^{\infty} )$ and often write $\pi_0(a)$ simply as $a$.

  We define 
 $ \Omega ( u ) := {\rm ev}_{e (u)} \circ \pi_\clq$, $\widetilde{\Omega}(u):={\rm ev}_{\tilde{e}(u)} \circ \pi_{\tilde{\clq}}$, for $u \in R^n$, where  ${\rm ev}_x$ (respectively ${\rm ev}_{\tilde{x}}$ )  denotes  the state on $C(\IT^n)$ (respectively, on $C(\tilde{\IT}^n)$) obtained by evaluation of a function at the point $x$ (respectively $\tilde{x}$). 
   
    
    For a fixed $J$, we shall work with several multiplications on the vector space ${\cla_0} \ot_{\rm alg} \clq_0$. We shall denote the counit and antipode of $\clq_0$ by $\epsilon$ and $\kappa$ respectively. Let us define the following

    $$ x \odot y = \int_{\IR^{4n}} e( -u.v )e( w.s )(\Omega ( -Ju )\diamond x \diamond (\Omega ( Jw ) ) (\Omega ( -v ) \diamond y \diamond \Omega( s )) du dv dw ds ,$$
     where $x,y \in \clq_0$. This is clearly a bilinear map, and will be seen to be an associative multiplication later on. 
    Moreover, we define two bilinear maps $\bullet$ and $\bullet_J$ by setting  $(a \ot x) \bullet (b \ot y):=ab \ot x \odot y$ and $(a \ot x) \bullet_J (b \ot y):=(a \times_J b) \ot (x \odot y)$, for $a,b \in {\cla_0}$, $x,y \in \clq_0$. 
We have $ \Omega(u) \diamond   ( \Omega(v) \diamond c ) = ( \Omega(u) \diamond \Omega(v) ) \diamond c  $. 

 
   \blmma 
    \label{Lemma1}

 1. The map $\odot$  satisfies $$ \int_{\IR^{2n}} ( \Omega( J u ) \diamond x ) \odot ( \Omega( v )   \diamond y )e(u.v ) du dv = \int_{\IR^{2n}} (x \diamond  ( \Omega( J u ))( y \diamond \Omega( v ) )e( u.v ) du dv ,$$
   for $x,y \in \clq_0$.
   
 2. $$ \widetilde{\alpha} ( \widetilde{\beta} _u ( s ) ) = s_{(1)} \otimes ( id \otimes \Omega( u ) ) ( \widetilde{\Delta} ( s_{( 2 )} ).$$
   
   $$ \alpha ( \beta_{u}( a ) ) = a_{1} \otimes ( id \otimes \Omega( u ) ) \widetilde{\Delta} ( a_{2} ) . $$
   
 3. For $ s \in {\cls}, a \in \widetilde{Q_{0}}$, we have      
     $$ \widetilde{\alpha} ( a \times_{J} s ) = a_{(1)}s_{(1)} \otimes ( \int\int ( a_{(2)} \diamond  Ju )  ( s_{(2)} \diamond v  ) e ( u.v ) du dv  ).$$  
   
 4.  For $ s \in {\cls_0}, a \in {\cla_{0}} $,       
     $$ \alpha ( a ) \bullet_J \widetilde{\alpha} ( s ) = a_{(1)}s_{(1)} \otimes \{ \int\int ( \Omega( Ju )  \diamond a_{(2)} ) \odot (  \Omega( v )  \diamond s_{(2)}) e( u.v ) du dv \}.$$ 
     
 5. For $a \in {\cla_0},~ s \in S $ we have $\alpha(a) \bullet_J \alpha(s)=\widetilde{\alpha}(a \times_J s).$

   \elmma
    {\it Proof :}\\ The proofs follow verbatim those in Lemmas 3.2 - 3.6 respectively in \cite{jyotish}. \qed

   \vspace{4mm}

 Let us recall at this point the Rieffel-type deformation of compact quantum groups as in \cite{wang_def}. We shall now identify $\odot$ with the multiplication of a Rieffel-type deformation of $\clq$. Since $\clq$ has a quantum subgroup isomorphic with $\IT^n$, we can consider the following canonical action $\lambda$ of $\IR^{2n}$ on $\clq$ given by          
      $$  \lambda_{( s,u )} = ( \Omega( -s ) \otimes id ) \Delta ( id \otimes \Omega ( u ) ) \Delta.$$  
          Now, let 
       $ \widetilde J := -J \oplus J $, which is a skew-symmetric ${2n} \times 2n$ real matrix, so one can deform $\clq$ by defining the product of $x$ and $y$ ($x,y \in \clq_0$, say) to be the following: $$ \int\int \lambda_{\widetilde J ( u,w )}( x ) \lambda _{v,s} ( y ) e ( ( u,w ).( v,s ) ) d ( u,w ) d ( v,s ).$$ We claim that this is nothing but $\odot$ introduced before.
         
         \blmma
              
           $ x \odot y = x \times_{\widetilde{J}} y ~ \forall x,y \in Q_{0} $
           
           \elmma
           
       {\it Proof :}\\ The proof is the same as Lemma 3.7 in \cite{jyotish}. \qed  
       
       \vspace{4mm}

     Let us denote by $\clq_{\widetilde{J}}$ the $C^*$ algebra obtained from $\clq$ by the  Rieffel deformation w.r.t. the matrix $\widetilde{J}$ described above. It has been shown in \cite{wang_def} that the coproduct $\Delta$ on $\clq_0$ extends to a coproduct for the deformed algebra  as well and  $(\clq_{\widetilde{J}}, \Delta)$ is a compact quantum group.

    We recall Lemma 3.8  of \cite{jyotish}, which is stated below for reader's convenience:  
     
      \blmma           
         \label{Lemma5a}
         
     The Haar state (say $h$) of $\clq$ coincides with the Haar state on $\clq_{\widetilde{J}}$ ( say $ h_{J} $ ) on the common subspace $\clq^{\infty}$, and moreover, $h(a \times_{\widetilde{J}} b)=h(ab)$ for $a,b \in \clq^{\infty}$.
         \elmma
         We note a useful implication of the above lemma. Let us make use of the identification of $\clq_0$ as a common vector-subspace of all $\clq_{\widetilde{J}}$. To be precise, we shall sometimes denote this identification map from $\clq_0 $ to $\clq_{\widetilde{J}}$ by $\rho_J$. 
\bcrlre
\label{repdef}
Let $W$ be a finite-dimensional (say, $n$-dimensional) unitary representation of $\clq$, with $\tilde{W} \in M_ n(\IC) \ot \clq_0$ be the corresponding unitary. Then, for any $\tilde{J}$, we have that  $\tilde{W}_J:=({\rm id} \ot \rho_J)(\tilde{W}) $ is unitary in $\clq_{\widetilde{J}}$, giving a unitary $n$-dimensional representation of $\clq_{\widetilde{J}}$.   In other words, any finite dimensional unitary representation of $\clq$ is also a unitary representation of $\clq_{\widetilde{J}}$.
\ecrlre
     {\it Proof:}\\
 Since the coalgebra structures of $\clq$ and $\clq_{\widetilde{J}}$ are identical, and $\tilde{W}_J$ is identical with $\tilde{W}$ as a linear map, it is obvious that $\tilde{W}_J$ gives a nondegenerate representation of $\clq_{\widetilde{J}}$. Let $y=({\rm id} \ot h)(\tilde{W}_J^* \tilde{W}_J)$. 
It follows from the proof of Proposition 6.4 of \cite{vandaelenotes} that $y$ is invertible positive element of $M_n$ and $(y^{\frac{1}{2}} \ot 1) \tilde{W}_J (y^{-\frac{1}{2}} \ot 1)$ gives a unitary representation of $\clq_{\widetilde{J}}$. We claim that $y=1$, which will complete the proof of the corollary.
For convenience, let us write $W$ in the Sweedler notation: $W=w_{(1)} \ot w_{(2)}$. We note that by Lemma \ref{Lemma5a}, we have \bean \lefteqn{({\rm id} \ot h)(\tilde{W}_J^*\tilde{W}_J)}\\
&=& w_{(1)}^*w_{(1)} h(w_{(2)}^* \times_{\widetilde{J}} w_{(2)})\\
&=& w_{(1)}^*w_{(1)} h(w_{(2)}^* w_{(2)})\\
&=& ({\rm id} \ot h)(\tilde{W}^* \tilde{W})=({\rm id} \ot h)(1 \ot 1)=1.\eean
\qed

\vspace{4mm}

Let us consider the finite dimensional unitary representations $U^{(i)}:=U|_{V_i}$, where $V_i$ is the eigenspace of $D$ corresponding to the eigenvalue $\lambda_i$. By the above Corollary \ref{repdef}, we can view $U^{(i)}$ as a unitary representation of $\clq_{\widetilde{J}}$ as well, and let us denote it by $U^{(i)}_J$. In this way, we obtain a unitary representation $U_J$ on the Hilbert space $\clh$, which is the closed linear span of all the $V_i$'s. It is obvious from the construction (and the fact that the linear span of $V_i$'s, i.e. $\cls_0$, is a core for $D$) that $U_J D=(D \ot I)U_J$.  Let $\alpha_J:=\alpha_{U_J}$.  
With this, we have the following:

  \blmma
For $ a \in \cla_0$, we have $\alpha_J(a)=(\alpha(a))_J \equiv (\pi_J \ot \rho_J)(\alpha(a))$, hence in particular, for every state $\phi$ on $\clq_{\widetilde{J}}$, $({\rm id} \ot \phi) \circ \alpha_J(\cla_J) \subseteq \cla_J^{\prime \prime}$.

  \elmma  
  
  Using Lemma \ref{Lemma1}, we have, $ \forall s \in \cls_0, a \in \cla_0$,
   \bean
  \lefteqn{ {U_J} ( \pi_J(a)s) } \\
       &=& \widetilde{\alpha}(a \times_J s)\\
&=& \alpha(a) \bullet_{J} \widetilde{\alpha}(s)\\
&=& ( \alpha(a) )_J U_J(s),
\eean
 from which we conclude by the density of $\cls_0$ in $\clh$ that $\alpha_J(a) =(\alpha(a))_J \in \pi_J(\cla_0) \ot \clq_{\widetilde{J}}$. The lemma now follows  using the norm-density of $\cla_0$ in $\cla_J$. 
\qed 

Thus, $(\tilde{\clq}_{\widetilde{J}}, U_J)$ is an orientation preserving isometric action on the spectral triple $(\cla^\infty_J, \clh, D)$.
    
    We shall now show that if we fix a `volume-form' in ternms of an $R$-twisted structure, then the `deformed' action $\alpha_J$ preserves it.
    
  \blmma
       \label{volumepreserving}
    Suppose, in adition to the set-up already assumed, that  there is an invertible positive operator $R$ on $\clh$ such that $(\cla^\infty, \clh, D, R)$ is an $R$-twisted $ \Theta $ summable spectral triple, and let $\tau_R$ be the corresponding `volume form'. Assume that  $\alpha_U$  preserves 
 the functional $\tau_R$. Then the action $\alpha_{U_J}$  preserves $\tau_R$ too.

    \elmma
      
  {\it Proof :}\\
             Let the (finite dimensional) eigenspace corresponding to the eigenvalue $ \lambda_n $ of $D$ be $ V_n .$   As $ U $ commutes with $ D ,$ there exists subspaces 
$ V_{n,k} $ of $ V_n $ and an orthonormal basis $ {\{e^{n,k}_{j}\}}_{j} $ for $ V_{n,k} $  such that the restriction of $U$ to $V_{n,k}$  is irreducible. Write $\widetilde{U}( e^{n,k}_{j} \otimes 1 ) = \sum_{i} e^{n,k}_{i} \otimes t^{n}_{i,j}.$ Then, $ {\widetilde{U}}^{*} ( e^{n,k}_{j} ) = \sum_{n,i} e^{n,k}_{i} \otimes t^{n *}_{j,i} .$
             
 Then $ \clh $ will be decomposed as $ \clh = \oplus_{n \geq 1, ~ k} V_{n,k} .$
             
     Let $ R ( e^{n,i}_{j} ) = \sum_{s,t} F_n ( i,j,s,t ) e^{n,s}_{t} .$
     
    
 By hypothesis, $ \widetilde{U} ( . \otimes id ) {\widetilde{U}}^{*} $ preserves the functional $ \tau_{R} ( \cdot) = Tr ( R~ \cdot ) $ on $ \cle_{D} $ where $ \cle_{D} $ is as in Proposition \ref{5678}, i.e, the weakly dense $ \ast $ subalgebra of $ \clb ( \clh ) $ generated by the rank one operators $ | \xi > < \eta | $ where $ \xi, \eta $ are eigenvectors of $ D .$      
   Thus, $( \tau_{R} \otimes id )( \widetilde{U} ( X \otimes id ){\widetilde{U}}^{\ast} ) = \tau_{R}( X ).1_{Q} \forall X \in \cle_{D}.$

   Then, for $a \in \cle_D$, we have: \bean
         \lefteqn{ ( \tau_{R} \otimes h )( \widetilde{U}_{J} ( a \otimes 1 ) {\widetilde{U}_{J}}^{*} ) } \\ 
         &=& \sum_{n,i,j} \left\langle e^{n,i}_{j} \otimes 1,\widetilde{U}_{J} ( a \otimes 1 ) {\widetilde{U}_{J}}^{*} ( R   e^{n,i}_{j} \otimes 1 ) \right\rangle  \\
         & = &  \sum_{n,i,j,s,t}  \left\langle  {\widetilde{U}_{J}}^{*} ( e^{n,i}_{j} \otimes 1 ), ( a \otimes 1 ){\widetilde{U}_{J}}^{*} (  F_{n}( i,j,s,t ) e^{n,s}_{t} \otimes 1 ) \right\rangle \\
         & = & \sum_{n,i,j,s,t,k,l} F_{n}(i,j,s,t ) \left\langle e^{n,i}_{k} \otimes {(t^{n}_{j,k})}^{*}, ( a \otimes 1 )( e^{n,s}_{l}  \otimes {(t^{n}_{t,l})}^{*} \right\rangle  \\
          &=& \sum_{n,i,j,s,t,k,l} F_{n}( i,j,s,t ) \left\langle e^{n,i}_{k} , a e^{n,s}_{l} \right\rangle h_{J} ( ( t^{n}_{j,k}) \times_{J} {(t^{n}_{t,l})}^{*}) \\
        & = &  \sum_{n,i,j,s,t,k,l} F_{n}( i,j,s,t ) \left\langle e^{n,i}_{k} , a e^{n,s}_{l} \right\rangle h_{0} ( t^{n}_{j,k} {t^{n}_{t,l}}^{*})  \\
   &=& ( \tau_{R} \otimes h )( \widetilde{U} ( a \otimes 1 ){ \widetilde{U} }^{*} ) \\ 
   &=&  \tau_{R} ( a ).1 \eean 
  
  where $ h_{J} ( ( t^{n}_{j,k}) \times_{J} {(t^{n}_{t,l})}^{*}) = h_{0} ( t^{n}_{j,k} {t^{n}_{t,l}}^{*}) $ as  deduced by using Lemma \ref{Lemma5a}.
  
  Thus $ ( \tau_R \otimes h ) \widetilde{U_{J}} ( a \otimes id ) { \widetilde{U_{J}} }^{*}  = \tau_R ( a ).1 $
  
  Let $ ( \tau_R \otimes h ) \widetilde{U_{J}} ( X \otimes id ) { \widetilde{U_{J}} }^{*}  = ( \tau_R \ast h ) ( X ) .$  As $ \widetilde{U_{J}} ( \cdot  \otimes id ) { \widetilde{U_{J}} }^{*} $ keeps $ \cle_{D} $ invariant, we can use Sweedler notation: $\widetilde{U_J}(a \ot 1) \widetilde{U_J}^*=
 a_{(1)} \ot a_{(2)}$, with $a, a_{(1)} \in \cle_D,~a_{(2)} \in \widetilde{\clq_{\tilde{J}}}$,   to have 
\bean \lefteqn{ ( \tau_{R} \otimes id ) (\widetilde{U_{J}} ( a \otimes 1 ) { \widetilde{U_{J}} }^{*})}\\
& =&
 ( \tau_{R} \ast h \otimes {\rm  id} )(\widetilde{U_{J}} ( a \otimes 1 ) { \widetilde{U_{J}} }^{*})  =  \tau_{R} \ast h ( a_{(1)})  a_{(2)}\\
& =&  ( \tau_{R} \otimes h  \ot {\rm id})  ( a_{(1)_{(1)}} \ot a_{(1)_{(2)}} \otimes a_{(2)}) = 
 ( \tau_{R} \otimes h  \otimes {\rm id} )( {\rm id } \otimes \Delta_{\widetilde{J}} ) (\widetilde{U_J}( a \ot 1 )\widetilde{U_J}^*)\\
 &= & \tau_{R}( a_{(1)} )  ( h \otimes {\rm  id} )\circ \Delta_{\widetilde{J}} ( a_{(2)} )
 = \tau_{R}( a_{(1)} )  h( a_{(2)} ).1_{\clq_{\widetilde{J}}} \\
&=&  ( \tau_{R} \otimes h ) ( a_{(1)} \ot a_{(2)} ) = ( \tau_{R} \ast h )( a ).1_{\clq_{\widetilde{J}}} = \tau_{R}( a ).1_{\clq_{\widetilde{J}}}.\eean

   \qed 
 
 \brmrk
   
 If  $ QISO^{+}_{R}( \cla^{\infty}, \clh, D ) $ (  $ QISO^{+}( \cla^{\infty}, \clh, D ),$ if it exists )  has a $ C^* $ action, then from the definition of a $ C^{*} $ action, we get a subalgebra $ \cla_0  $ as in Lemma \ref{1111}. Thus, the conclusions of Lemma \ref{1111} and the subsequent Lemmas hold for $ QISO^{+}_{R}( \cla^{\infty}, \clh, D ) $ (  $ QISO^{+}( \cla^{\infty}, \clh, D ) .$
 
 \ermrk

 
We have already seen that $(\tilde{\clq}_{\widetilde{J}}, U_J)$ is an object in ${\bf Q} ( \cla_J, \clh, D )$.
 Now, proceeding as in the proof of Theorem 3.13 of \cite{jyotish} we obtain the following result ( using Lemma  \ref{volumepreserving} for 1 ). 
     
\bthm
\label{main_def}

1. If $ QISO^{+}_{R}(\cla^\infty_J, \clh, D) $ and $ (QISO^{+}_{R}(\cla^\infty, \clh, D))_{\widetilde{J}} $ have $ C^* $ actions on $ \cla $ and $ \cla_J $ respectively, we have

$$ \widetilde{QISO^{+}_{R}} ( \cla^{\infty}_{J}, \clh, D ) \cong {( \widetilde{QISO^{+}_{R}} ( \cla^{\infty}, \clh, D ) )}_{\widetilde{J}}  $$

$$ QISO^{+}_{R}(\cla^\infty_J, \clh, D) \cong (QISO^{+}_{R}(\cla^\infty, \clh, D))_{\widetilde{J}}.$$

2. If moreover, $\widetilde{QISO}^+(\cla^\infty, \clh, D)$ and $ \widetilde{{QISO}^{+}} ( \cla^{\infty}_{J}, \clh, D ) $ both exist and have $ C^* $ actions on $ \cla $ and  $ \cla_{J} $ respectively, then 
  $${\widetilde{QISO}}^+(\cla^\infty_J, \clh, D) \cong \left({\widetilde{QISO}}^+(\cla^\infty, \clh, D)\right)_{\widetilde{J}},$$  $$    QISO^+(\cla^\infty_J, \clh, D) \cong (QISO^+(\cla^\infty, \clh, D))_{\widetilde{J}}.$$

\ethm

\vspace{8mm}

 As an example, we consider the noncommutative torus $ \cla_{\theta} ,$ which is  a Rieffel deformation of  $ C( \IT^{2} ) $ with respect to the matrix $ J =   \left ( \begin {array} {cccc}
    0 &  \theta   \\ - \theta &  0  \end {array} \right ) $
    and we deform the spectral triple as in subsection 4.This is the standard spectral triple on $ \cla_{\theta}.$
   
   \bthm
    
     $\widetilde{QISO}^+(\cla^\infty_\theta, \clh, D ) = \widetilde{QISO}^+(C^\infty(\IT^2))=C(\IT^2) \ast C(\IT)$, and $QISO^+(\cla^\infty_\theta)=QISO^+(C^\infty(\IT^2))=C(\IT^2)$.
     
     \ethm
    
{\it Proof:}\\   We use Theorem \ref{main_def} and recall that $ QISO^+(C^\infty(\IT^2))=C(\IT^2)$ which is generated by $ z_{1} $ and $ z_{2} ,$ say.

Then, from the formula of the deformed product, it can easily be seen after a change of variable that $ z_{1}  \times_{\widetilde{J}} ~  z_{2} = z_{2} \times_{\widetilde{J}} ~ z_{1} $ which proves the theorem.

\qed

\brmrk
In a private communication  S. Wang has kindly pointed out that one can possibly formulate and prove an analogue of Theorem \ref{main_def}
 in the setting of discrete deformation as in \cite{wang_disc}, and this may give a solution to a problem posed by Connes (see \cite{con}, page 612). We believe
 that more work is needed in this direction.
\ermrk

{\bf Acknowledgement:}\\
We would like to thank S. Wang and T. Banica for many valuable comments and feedback which have led to substantial improvement of an earlier version of this article.  \vspace{2mm}\\

\end{document}